\documentclass[12pt]{article}

\usepackage{amssymb}
\usepackage{amsmath}
\usepackage{graphicx}
\usepackage{pst-all}
\usepackage{mathrsfs}
\usepackage{cite}
\usepackage{pst-plot}  
\usepackage{tikz}
\usepackage{pst-poly}

\newtheorem{thm}{Theorem}[section]
\newtheorem{cor}[thm]{Corollary}
\newtheorem{lem}[thm]{Lemma}

\newenvironment{pf}[1][Proof]{\noindent\textbf{#1.} }{\hfill\rule{1mm}{2mm}}

\begin{document}
\title
{Relationship between Conditional Diagnosability and 2-extra
Connectivity of Symmetric Graphs\thanks{The work was supported by
NNSF of China No. 11371052, 11171020, 11271012, 61272008.} }

\author
{Rong-Xia Hao$^a$\quad Zeng-Xian Tian$^a$ and Jun-Ming
Xu$^b$\footnote{Corresponding author: xujm@ustc.edu.cn (J.-M. Xu)}
\footnote{Email addresses: rxhao@bjtu.edu.cn (R.-X. Hao), 14121550@bjtu.edu.cn (Z.-X. Tian).}\\
{\small $^a$Department of Mathematics, Beijing Jiaotong University, Beijing, 100044, China}\\
\small $^b$School of Mathematical Sciences, University of Science and Technology of China, \\
\small Hefei, Anhui, 230026, China\\
}

\date{} \maketitle

\begin{minipage}{140mm}

\begin{center} {\bf Abstract} \end{center}

The conditional diagnosability and the 2-extra connectivity are two
important parameters to measure ability of diagnosing faulty
processors and fault-tolerance in a multiprocessor system. The {\it
conditional diagnosability} $t_c(G)$ of $G$ is the maximum number
$t$ for which $G$ is conditionally $t$-diagnosable under the
comparison model, while the {\it $2$-extra connectivity}
$\kappa_2(G)$ of a graph $G$ is the minimum number $k$ for which
there is a vertex-cut $F$ with $|F|=k$ such that every component of
$G-F$ has at least $3$ vertices. A quite natural problem is what is
the relationship between the maximum and the minimum problem? This
paper partially answer this problem by proving $t_c(G)=\kappa_2(G)$
for a regular graph $G$ with some acceptable conditions. As
applications, the conditional diagnosability and the 2-extra
connectivity are determined for some well-known classes of
vertex-transitive graphs, including, star graphs, $(n,k)$-star
graphs, alternating group networks, $(n,k)$-arrangement graphs,
alternating group graphs, Cayley graphs obtained from transposition
generating trees, bubble-sort graphs, $k$-ary $n$-cube networks and
dual-cubes. Furthermore, many known results about these networks are
obtained directly.

\vskip0.4cm \noindent

{\bf Keywords}\ conditional diagnosability; comparison model; extra
connectivity; symmetric graph; Cayley graph; max-min problem


\end{minipage}

\setlength{\baselineskip}{18pt}

 \vskip1.2cm

\section{Introduction}

Throughout this paper, unless otherwise specified, a graph $G=(V,E)$
is always assumed to be a simple and connected graph, where $V=V(G)$
is the vertex-set and $E=E(G)$ is the edge-set of $G$. We
follow~\cite{x03} for terminologies and notations not defined here.

Two distinct vertices $x$ and $y$ in $G$ is adjacent if $xy\in E(G)$
and non-adjacent otherwise. If $xy\in E(G)$, then $y$ (resp. $x$) is
a neighbor of $x$ (resp. $y$). The neighbor-set of $x$ is denoted by
$N_G(x)=\{y\in V(G):\ xy\in E(G)\}$. For a subset $X\subset V(G)$,
the notation $G-X$ denotes the subgraph obtained from $G$ by
deleting all vertices in $X$ and all edges incident with vertices in
$X$, and let $\overline{X}=V(G-X)$.

It is well known that a topological structure of an interconnection
network $N$ can be modeled by a graph $G=(V,E)$, where $V$
represents the set of components such as processors and $E$
represents the set of communication links in $N$ (see a text-book by
Xu~\cite{x13}). Faults of some processors and/or communication lines
in a large-scale system are inevitable. People are concerned with
how to diagnose faults and to determine fault tolerance of the
system.

A vertex in a graph $G$ is called a {\it fault-vertex} if it
corresponds a faulty processor in the interconnection network $N$
when it is modeled by $G$. A subset $F\subseteq V(G)$ is called a
{\it fault-set} if every vertex in $F$ is a faulty vertex in $G$,
and is {\it fault-free} if it contains no faulty vertex in $G$. A
fault-set $F$ is called a {\it conditional fault-set} if
$N_G(x)\nsubseteq F$ for any $x\in \overline{F}$. The pair
$(F_1,F_2)$ is called a {\it conditional fault-pair} if both $F_1$
and $F_2$ are conditional fault-sets.

The ability to identify all faulty processors in a multiprocessor
system is known as system-level diagnosis. Several system-level
self-diagnosis models have been proposed for a long time. One of the
most important models is the {\it comparison diagnosis model},
shortly {\it comparison model}. Throughout this paper, we only
consider the comparison model.

The comparison model was proposed by Malek and Maeng~\cite{mm81,
m80}. A node can send a message to any two of its neighbors which
then send replies back to the node. On receipt of these two replies,
the node compares them and proclaims that at least one of the two
neighbors is faulty if the replies are different or that both
neighbors are fault-free if the replies are identical. However, if
the node itself is faulty then no reliance can be placed on this
proclamation. According as that the two outputs are identical or
different, one gets the outcome to 0 or 1. The collection of all
comparison results forms a syndrome, denoted by $\sigma$.

A subset $F\subseteq V(G)$ is a {\it compatible\ fault-set} of a
syndrome $\sigma$ or $\sigma$ is {\it compatible with} $F$, if
$\sigma$ can arise from the circumstance that $F$ is a fault-set and
$\overline{F}$ is fault-free. Let $\sigma_F=\{{\sigma:\ \sigma\ is \
\mbox{compatible\ with}\ F}\}$. A pair $(F_1,F_2)$ of two distinct
compatible fault-sets is {\it distinguishable} if and only if
$\sigma_{F_1}\cap\sigma_{F_2}=\emptyset$, and $(F_1, F_2)$ is {\it
indistinguishable} otherwise.

For a positive integer $t$, a graph $G$ is {\it conditionally
$t$-diagnosable} if every syndrome $\sigma$ has a unique conditional
compatible fault-set $F$ with $|F|\leqslant t$. The {\it conditional
diagnosability} of $G$ under the comparison model, denoted by
$t_c(G)$ and proposed by Lai et al. {\rm\cite{ltch05}}, is the
maximum number $t$ for which $G$ is conditionally $t$-diagnosable.
The conditional diagnosability better reflects the self-diagnostic
capability of networks under more practical assumptions, and has
received much attention in recent years. The diagnosability of many
interconnection networks have been determined, see, for example,
\cite{as00, as02, f98a, f98b, f02, hfz13, kk91, w94}. A survey on
this field, from the earliest theoretical models to new promising
applications, is referred to Duarte {\it et al.}~\cite{dza11}.

A subset $X\subset V(G)$ is called a {\it vertex-cut} if $G-X$ is
disconnected. A vertex-cut $X$ is called a $k$-cut if $|X|=k$. The
{\it connectivity} $\kappa(G)$ of $G$ is defined as the minimum
number $k$ for which $G$ has a $k$-cut.

Fault-tolerance or reliability of a large-scale parallel system is
often measured by the connectivity $\kappa(G)$ of a corresponding
graph $G$. However, the connectivity has an obvious deficiency
because it tacitly assumes that all vertices adjacent to the same
vertex of $G$ could fail at the same time, but that is almost
impossible in practical network applications. To compensate for this
shortcoming, F\`{a}brega and Fiol {\rm\cite{ff96}} proposed the
concept of the extra connectivity.

For a non-negative positive integer $h$, a vertex-cut $X$ is called
an {\it $R_h$-vertex-cut} if every component of $G-X$ has at least
$h+1$ vertices. For an arbitrary graph $G$, $R_h$-vertex-cuts do not
always exist for some $h$. For example, a cycle of order $5$
contains no $R_2$-vertex-cut. A graph $G$ is called an $R_h$-graph
if it contains at least one $R_h$-vertex-cut. For an $R_h$-graph
$G$, the {\it $h$-extra connectivity} of $G$, denoted by
$\kappa_h(G)$, is defined as the minimum number $k$ for which $G$
contains an $R_h$-vertex-cut $F$ with $|F|=k$. Clearly,
$\kappa_0(G)=\kappa(G)$. Thus, the $h$-extra connectivity is a
generalization of the classical connectivity and can provide more
accurate measures regarding the fault-tolerance or reliability of a
large-scale parallel system and therefore, it has received much
attention (see Xu~\cite{x13} for details). We are interested in the
$2$-extra connectivity of a graph in this paper.

Clearly, for a graph $G$ there are two problems here, one is the
maximizing problem -- conditional diagnosability $t_c(G)$, and
another is the minimizing problem -- the $2$-extra connectivity
$\kappa_2(G)$. A quite natural problem is what is the relationship
between the maximum and the minimum problems? In the current
literature, people are still determining these two problems
independently for some classes of graphs, such as alternating group
network \cite{z09}, alternating group graph \cite{hz12, zxy10,
zx10}, the $3$-ary $n$-cube network \cite{zj13}.

In this paper, we reveal the relationships between the conditional
diagnosability $t_c(G)$ and the $2$-extra connectivity $\kappa_2(G)$
of a regular graph $G$ with some acceptable conditions by
establishing $t_c(G)=\kappa_2(G)$. As applications of our result, we
consider some more general well-known classes of vertex-transitive
graphs, such as star graphs, $(n,k)$-star graphs, alternating group
networks, $(n,k)$-arrangement graphs, alternating group graphs,
Cayley graphs obtained from transposition generating trees,
bubble-sort graphs, $k$-ary $n$-cube networks and dual-cubes, and
obtain the conditional diagnosability under the comparison model and
the $2$-extra connectivity of these graphs, which contain all known
results on these graphs.

The rest of the paper is organized as follows. Section 2 first
recalls some necessary notations and lemmas, then establishes the
relationship between the conditional diagnosability and the 2-extra
connectivity of regular graphs with some conditions. As applications
of our main result, Section 3 determines the conditional
diagnosability and the 2-extra connectivity for some well-known
classes of vertex-transitive graphs.

\section{Main results}

We first recall some terminologies and notation used in this paper.
Let $G=(V,E)$ be a graph, where $V=V(G)$, $E=E(G)$ and $|V(G)|$ is
the order of $G$.

A sequence $(x_1, \cdots, x_n)$ of $n\ (\geqslant 3)$ distinct
vertices with $x_ix_{i+1}\in E(G)$ for each $i=1,\cdots,n-1$ is
called an $n$-path, denoted by $P_n$, if $x_1x_n\notin E(G)$, and
called an $n$-cycle, denoted by $C_n$, if $x_1x_n\in E(G)$. A cycle
$C$ in $G$ is {\it chordless} if any two non-adjacent vertices of
$C$ are non-adjacent in $G$.

For $X\subset V(G)$, let $N_G(X)=(\cup_{x\in X}\ N_G(x))\setminus
X$. For simplicity of writing, in case of no confusion from the
context, we write $N(x)$ for $N_G(x)$; moreover, if $X$ is a
subgraph of $G$, we write $N(X)$ for $N_G(V(X))$ in this paper. For
two non-adjacent vertices $x$ and $y$ in $G$, let
$\ell(x,y)=|N(x)\cap N(y)|$, and let $\ell(G)=\max\{\ell(x,y):\
x,y\in V(G)\ {\rm and}\ xy\notin E(G)\}$.

The degree $d(x)$ of a vertex $x$ is the number of neighbors of $x$,
i.e., $d(x)=|N(x)|$. The minimum degree $\delta(G)=\min\{d(x):\ x\in
V(G)\}$ and the maximum degree $\Delta(G)=\max\{d(x):\ x\in V(G)\}$.
A vertex $x$ is an {\it isolated vertex} if $d(x)=0$, an edge $xy$
is an {\it isolated edge} if $d(x)=d(y)=1$. A graph $G$ is
$k$-regular if $\delta(G)=\Delta(G)=k$. $K_n$ denotes a complete
graph of order $n$, which is an $(n-1)$-regular graph. For a
subgraph $H$ of $G$, we will use $\Sigma(H)$ to denote $\Sigma_{x\in
H} d_H(x)$. For example, if $P_3$ and $C_3$ are subgraphs of $G$,
then $\Sigma(P_3)=4$ and  $\Sigma(C_3)=6$.

Let $X\subset V(G)$ be a vertex-cut. The maximal connected subgraphs
of $G-X$ are called {\it components}. A component is {\it small} if
it is an isolated vertex or an isolated edge; is {\it large}
otherwise.


In this section, we present our main theorem, which explores the
close relationship between the conditional diagnosability $t_c(G)$
and the 2-extra connectivity $\kappa_2(G)$ of a regular graph $G$
under some conditions, that is, $t_c(G)=\kappa_2(G)$. The following
three lemmas play a key role in the proof of our theorem.

\begin{lem}{\rm\cite{sd92}}\label{lem2.1}
Let $G=(V,E)$ be a graph, $F_1, F_2\subseteq V(G), F_1\neq F_2$.
Then, under the comparison model, $(F_1,F_2)$ is a distinguishable
pair if and only if one of the following conditions is satisfied
{\rm (see Fig.~\ref{f1})}.

\vskip4pt

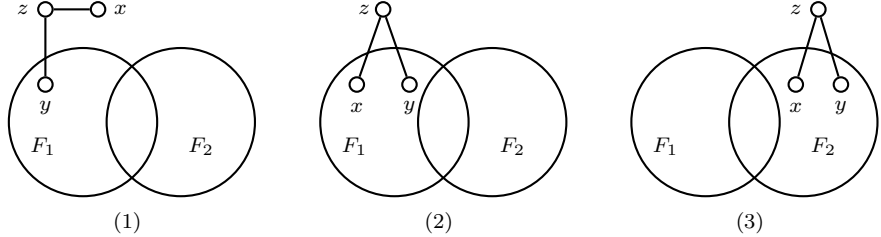
\begin{figure}[h]
\begin{pspicture}(-2.,-0.2)(1,2.5)

 \psellipse(1.,1)(1.,1.)
 \psellipse(2.3,1)(1.,1.)

 \cnode(.5,2.5){.11}{a1}\rput(.2,2.5){\scriptsize$z$}
 \cnode(.5,1.5){.11}{a2}\rput(.5,1.2){\scriptsize$y$}
 \cnode(1.2,2.5){.11}{a3}\rput(1.5,2.5){\scriptsize$x$}
 \ncline{a1}{a2}\ncline{a1}{a3}

 \put(0.3,.6){\scriptsize$F_1$}
 \put(2.4,0.6){\scriptsize$F_2$}
 \put(1.4,-0.4){\scriptsize$(1)$}

  \end{pspicture}
\begin{pspicture}(-3.,-0.2)(0,2.5)

 \psellipse(1.,1)(1.,1.)
 \psellipse(2.3,1)(1.,1.)

 \cnode(.85,2.5){.11}{a1}\rput(.6,2.5){\scriptsize$z$}
 \cnode(.5,1.5){.11}{a2}\rput(.5,1.2){\scriptsize$x$}
 \cnode(1.2,1.5){.11}{a3}\rput(1.2,1.2){\scriptsize$y$}
 \ncline{a1}{a2}\ncline{a1}{a3}

 \put(0.3,.6){\scriptsize$F_1$}
 \put(2.4,0.6){\scriptsize$F_2$}
 \put(1.4,-0.4){\scriptsize$(2)$}

  \end{pspicture}
\begin{pspicture}(-4.,-0.2)(0,2.5)

 \psellipse(1.,1)(1.,1.)
 \psellipse(2.3,1)(1.,1.)

 \cnode(2.5,2.5){.11}{b1}\rput(2.2,2.5){\scriptsize$z$}
 \cnode(2.2,1.5){.11}{b2}\rput(2.2,1.2){\scriptsize$x$}
 \cnode(2.8,1.5){.11}{b3}\rput(2.8,1.2){\scriptsize$y$}
 \ncline{b1}{b2}\ncline{b1}{b3}

 \put(0.3,.6){\scriptsize$F_1$}
 \put(2.4,0.6){\scriptsize$F_2$}
 \put(1.4,-0.4){\scriptsize$(3)$}

  \end{pspicture}

\caption{\label{f1} \footnotesize Illustrations of
Lemma~\ref{lem2.1}}
\end{figure}

{\rm (a)}\ There exists $x,z\in \overline{F_1\cup F_2}$ and $y\in
(F_1\cup F_2)\setminus (F_1\cap F_2)$ such that $xz,yz\in E(G)$;

{\rm (b)}\ There exists $z\in \overline{F_1\cup F_2}$ and $x,y\in
F_1\setminus F_2$ such that $xz,yz\in E(G)$;

{\rm (c)}\ There exists $z\in \overline{F_1\cup F_2}$ and $x,y\in
F_2\setminus F_1$ such that $xz,yz\in E(G)$.
\end{lem}

\begin{lem}{\rm\cite{sd92}}\label{lem2.2}
A graph $G$ is conditionally $t$-diagnosable if and only if, for any
two distinct conditional fault-sets $F_1$ and $F_2$ with
$\max\{|F_1|, |F_2|\}\leqslant t$, $(F_1,F_2)$ is a distinguishable
pair.
\end{lem}

\begin{lem}{\rm\cite{clqs12}}\label{lem2.3}\
Let $G=(V,E)$ be a graph with maximum degree $\Delta$ and minimum
degree $\delta\geqslant 3$. If there is some integer $t$ such that

{\rm (a)}\ $|V|>(\Delta+1)(t-1)+4$,{\footnote{This lower bound on
$|V|$ given here is quite enough for the conclusion. The original
article claims $|V|>(\Delta+2)(t-1)+4$. }}

{\rm (b)}\ for any $F\subset V(G)$ with $|F|\leqslant t-1$, $G-F$
has a large component and small components (if exist) which contain
at most two vertices in total.

\noindent then $t_c(G)\geqslant t$.
\end{lem}




\begin{thm}\label{thm2.4}
Let $G$ be an $n$-regular $R_2$-graph and\, $t=\min\{|N(T)|:\ T$ is
a $3$-path or a $3$-cycle in $G\}$. If $G$ satisfies the following
conditions

{\rm (a)}\ for any $F\subset V(G)$ with $|F|\leqslant t-1$, $G-F$
has a large component and small components which contain at most two
vertices in total,

{\rm (b)}\ $n\geqslant 2\ell(G)+2$ if $G$ contains no $5$-cycle, and
$n\geqslant 3\ell(G)+2$ otherwise,

{\rm (c)}\ $|V(G)|>(n+1)(t-1)+4$,

\noindent then $t_c(G)=t=\kappa_2(G)$.
 \end{thm}

\begin{pf}
Let $T=P_{3}$ or $C_3$ (if exists) in $G$ such that $|N(T)|=t$. The
condition (c) implies that $N(T)$ is a vertex-cut of $G$.

Suppose that $N(T)$ is not an $R_2$-vertex-cut of $G$. Then $G-N(T)$
contains a small component $C$ which contains at most two vertices.

If $C$ is an isolated vertex, say $x$, then $x$ shares at most
$\ell(G)$ common neighbors with any of three vertices in $T$. Thus,
$n=|N(x)\cap N(T)|\leqslant\min\{3\ell(G),n\}$, which implies
$n\leqslant 3\ell(G)$, a contradiction with the hypothesis (b) that
$n\geqslant 3\ell(G)+2$. Moreover, if $G$ contains no $5$-cycle,
then $x$ shares at most $\ell(G)$ common neighbors with each of at
most two vertices in $T$, and so $n=|N(x)\cap
N(T)|\leqslant\min\{2\ell(G),n\}$, which implies $n\leqslant
2\ell(G)$, a contradiction with the hypothesis (b) that $n\geqslant
2\ell(G)+2$.

If $C$ is an isolated edge, say $xy$, then at most $(n-1)$ neighbors
of $x$ are in $N(T)$. In the same discussion above, we have that
$n-1=|N(x)\cap N(T)|\leqslant\min\{3\ell(G),n-1\}$, which implies
$n\leqslant 3\ell(G)+1$; and if $G$ contains no $5$-cycle, then
$n-1=|N(x)\cap N(T)|\leqslant\min\{2\ell(G),n-1\}$, which implies
$n\leqslant 2\ell(G)+1$. These contradict with the condition (b).
Hence, $N(T)$ is an $R_2$-vertex-cut of $G$, and so
$\kappa_2(G)\leqslant |N(T)|=t$.

On the other hand, since $G$ is an $R_2$-graph, there is an
$R_2$-vertex-cut $F$ of $G$ such that $|F|=\kappa_2(G)$. Clearly,
$F$ is a vertex-cut of $G$. By the condition (a), if $|F|\leqslant
t-1$, then $G-F$ certainly contains a small component $C$ with
$|V(C)|\leqslant 2$, which contradicts the assumption that $F$ is an
$R_2$-vertex-cut, and so $\kappa_2(G)=|F|\geqslant t$. Thus,
$\kappa_2(G)=t$.

\vskip6pt

We now prove  $t_c(G)=t$. The conditions (a) and (c) satisfy two
conditions in Lemma~\ref{lem2.3}, and so $t_c(G)\geqslant t$.

On the other hand, let $T=\{x,z,y\}$ with $xz, yz\in E(G)$ such that
$|N(T)|=t$. By the above discussion, $N(T)$ is an $R_2$-vertex-cut
of $G$. Let $F_1=N(T)\cup \{x\}$ and $F_2=N(T)\cup \{y\}$. Then
$F_1\ne F_2$ and $|F_1|=|F_2|=t+1$. If there is a vertex $u\in
\overline{F_1}$ such that $N(u)\subseteq F_1$, then $u\notin\{y,z\}$
clearly, and so $u$ is in $G-N[T]$. Since $u$ is not adjacent to
$x$, $u$ is an isolated vertex in $G-N(T)$, which implies that
$N(T)$ is not an $R_2$-vertex-cut, a contradiction. Therefore, $F_1$
is a conditional fault-set. Similarly, $F_2$ is also a conditional
fault-set. Note that $(F_1\cup F_2)\setminus (F_1\cap F_2)=\{x,y\}$,
$F_1\setminus F_2=\{x\}$ and $F_2\setminus F_1=\{y\}$. It is easy to
verify that $F_1$ and $F_2$ satisfy none of conditions in Lemma
\ref{lem2.1}, and so $(F_1,F_2)$ is an indistinguishable pair. By
Lemma~\ref{lem2.2}, $G$ is not conditionally $(t+1)$-diagnosable,
which implies $t_c(G)\leqslant t$. Thus, $t_c(G)=t$.

It follows that $t_c(G)=t=\kappa_2(G)$. The theorem follows.
\end{pf}


\vskip20pt

\section{Applications to Some Well-known Networks}

As applications of Theorem \ref{thm2.4}, in this section, we
determine the conditional diagnosability and 2-extra connectivity
for some well-known vertex-transitive graphs, which, due to their
high symmetry, frequently appear in the literature on designs and
analyses of interconnection networks, including star graphs,
alternating group networks, alternating group graphs, bubble-sort
graphs, $(n,k)$-arrangement graphs, $(n,k)$-star graphs, a class of
Cayley graphs obtained from transposition generating trees, $k$-ary
$n$-cube networks and dual-cubes as well.

\subsection{Preliminary on Groups and Cayley Graphs}

We first simply recall some basic concepts on groups and the
definition of Cayley graphs, and introduce two classes of Cayley
graphs based on the alternating group, alternating group networks
and alternating group graphs.

Denote by $\Omega_n$ the group of all permutations on
$I_n=\{1,\ldots,n\}$. For convenience, we use $p_1p_2\cdots p_n$ to
denote the permutation ${1\ 2\ \cdots\ n}\choose{p_1 p_2 \cdots
p_n}$. A {\it transposition} is a permutation that exchanges two
elements and leaves the rest unaltered. A transposition that
exchanges $i$ and $j$ is denoted by $(i,j)$.

It is well known that any permutation can be expressed as
multiplications of a series of transpositions with operation
sequence from left to right.
In
particular, a 3-cycle $(a,b,c)$ is always expressed as
$(a,b,c)=(a,b)(a,c)$. For example, $(1,2,4)=(1,2)(1,4)$.

A permutation is called {\it even} if it can be expressed as a
composition of even transpositions, and {\it odd} otherwise. There
are $n!/2$ even permutations in $\Omega_n$, which form a subgroup of
$\Omega_n$, called the {\it alternating group} and denoted by
$\Gamma_n$, the generating set to be a set of 3-cycles.

\vskip6pt

An automorphism of a graph $G$ is a permutation on $V(G)$ that
preserves adjacency. All automorphisms of $G$ form a group, denoted
by Aut\,$(G)$, and referred to as the automorphism group. A graph
$G$ is {\it vertex-transitive} if for any two vertices $x$ and $y$
in $G$ there is a $\sigma \in$ Aut\,$(G)$ such that $y=\sigma (x)$.
A vertex-transitive graph is necessarily regular. A graph $G$ is
{\it edge-transitive} if for any two edges $a=xy$ and $b=uv$ of $G$
there is a $\sigma \in Aut(G)$ such that
$\{u,v\}=\{\sigma(x),\sigma(y)\}$. A graph is {\it symmetric} if it
is vertex-transitive and edge-transitive.

For a finite group $\Gamma$ with the identity $e$ and a non-empty
subset $S$ of $\Gamma$ such that $e\notin S$ and $S= S^{-1}$, define
a graph $G$ as follows.
$$
V(G)=\Gamma;\quad xy\in E(G)\Leftrightarrow x^{-1}y\in S\ {\rm for\
any}\ x,y\in \Gamma.
$$
In other words, $xy\in E(G)$ if and only if there exists $s\in S$
such that $y=xs$. Such a graph $G$ is called the {\em Cayley graph}
on $\Gamma$ with respect to $S$, denoted by $C_\Gamma(S)$. A Cayley
graph is $|S|$-regular, and is connected if and only if $S$
generates $\Gamma$. Moreover, A Cayley graph is $|S|$-connected if
$S$ is a minimal generating set of $\Gamma$.

A Cayley graph is always vertex-transitive and, thus, becomes an
important topological structure of interconnection networks and has
attracted considerable attention in the literature~\cite{hd97,
ljd93}.

\vskip6pt

As examples, we recall two well-known classes of Cayley graphs on
the alternating group $\Gamma_n$ with respect to some $S$.

\vskip10pt

{\bf 1. Alternating Group Networks}

\begin{figure}[h]
\begin{pspicture}(-2.5,-1.6)(1.7,3)
\psset{radius=.1}

\Cnode(.5,.5){231}\rput(0.2,.2){\scriptsize 231}
\Cnode(2.3,.5){312}\rput(2.5,.2){\scriptsize 312}
\Cnode(1.4,1.9){123}\rput(1.4,2.2){\scriptsize 123}
\ncline{231}{312}\ncline{231}{123}\ncline{123}{312}
\rput(1.25,-.2){\scriptsize $AN_3$}
\end{pspicture}
\begin{pspicture}(-1.5,0)(5,5.3)
\psset{radius=.1}

\Cnode(2.2,0.5){4321}\rput(2.,.2){\scriptsize4321}
\Cnode(5.,0.5){3412}\rput(5.2,0.2){\scriptsize3412}
\Cnode(1.4,1.6){2431}\rput(1.,1.6){\scriptsize2431}

\Cnode(5.8,1.6){4132}\rput(6.3,1.6){\scriptsize 4132}
\Cnode(2.7,1.6){3241}\rput(3.15,1.5){\scriptsize3241}
\Cnode(4.6,1.6){1342}\rput(4.15,1.5){\scriptsize1342}

\Cnode(2.8,4.9){4213}\rput(2.35,4.85){\scriptsize 4213}
\Cnode(4.3,4.9){1423}\rput(4.75,4.85){\scriptsize 1423}
\Cnode(3.6,4){2143}\rput(4.1, 3.9){\scriptsize 2143}

\Cnode(3.6,3.2){1234}\rput(4.1,3.2){\scriptsize 1234}
\Cnode(3.,2.3){2314}\rput(2.5,2.4){\scriptsize 2314}
\Cnode(4.2,2.3){3124}\rput(4.7,2.4){\scriptsize 3124}

\ncline{4321}{3241}\ncline{3241}{2431}\ncline{2431}{4321}
\ncline{3412}{4132}\ncline{4132}{1342}\ncline{1342}{3412}
\ncline{4213}{1423}\ncline{1423}{2143}\ncline{2143}{4213}
\ncline{1234}{2314}\ncline{2314}{3124}\ncline{3124}{1234}

\ncline{4321}{3412}\ncline{4132}{1423}\ncline{4213}{2431}
\ncline{1234}{2143}\ncline{2314}{3241}\ncline{1342}{3124}

\rput(3.6,0){\scriptsize $AN_4$}

\end{pspicture}
\caption{\label{f2}\footnotesize{Alternating group networks $AN_3$
and $AN_4$}}
\end{figure}
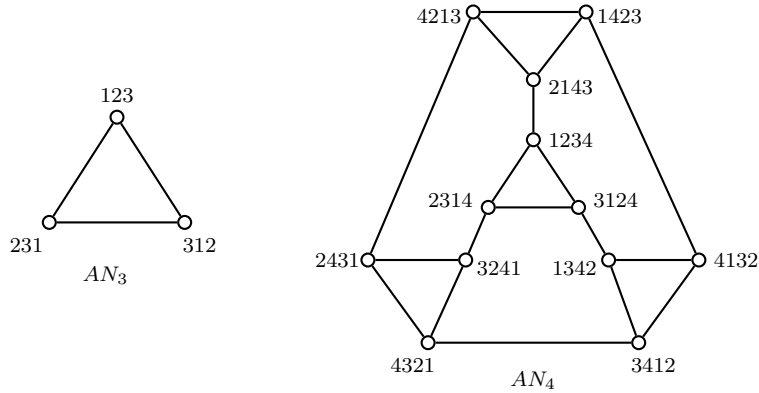

For $n\geqslant 3$, let $S=\{(1,2)(1,3), (1,3)(1,2),$ $(1,2)(3,i):$
$\ 4\leqslant i\leqslant n\}$, where $(1,2)(1,3)$ and $(1,3)(1,2)$
are mutually inverse, $(1,2)(3,i)$ is self-inverse for each
$i=4,\cdots,n$, and so $S=S^{-1}$. The Cayley graph
$C_{\Gamma_n}(S)$ is called the {\em alternating group network},
proposed by Ji~\cite{j99} in 1999 and denoted by $AN_n$, which is
$(n-1)$ regular and $(n-1)$-connected. The alternating group
networks $AN_{3}$ and $AN_{4}$ are shown in Fig.~\ref{f2}.

Zhou and Xiao~\cite{zx10} determined $t_c(AN_n)=3n-9$ for
$n\geqslant 5$ and Zhou~\cite{z09} determined $\kappa_2(AN_n)=3n-9$
for $n\geqslant 4$. Thus, $t_c(AN_n)=3n-9=\kappa_2(AN_n)$ for
$n\geqslant 5$

\vskip10pt

 {\bf 2. Alternating Group Graphs}

\begin{figure}[h]
\begin{pspicture}(-2,-1.6)(2.5,3)
\psset{radius=.1}

\Cnode(.5,.5){231}\rput(0.2,.2){\scriptsize 231}
\Cnode(2.3,.5){312}\rput(2.5,.2){\scriptsize 312}
\Cnode(1.4,1.9){123}\rput(1.4,2.2){\scriptsize 123}
\ncline{231}{312}\ncline{231}{123}\ncline{123}{312}

\rput(1.25,-.2){\scriptsize $AG_3$}
\end{pspicture}
\begin{pspicture}(-1.,0)(5,5.)
\psset{radius=.1}
 \Cnode(2.5,0.5){4213}\rput(2.1,.3){\scriptsize 4213}
 \Cnode(4.7,0.5){1423}\rput(5.1,0.3){\scriptsize 1423}
 \Cnode(3.6,1.2){2143}\rput(3.6,.85){\scriptsize 2143}
 \Cnode(2.5,1.9){2314}\rput(2.9,1.65){\scriptsize 2314}
 \Cnode(4.7,1.9){3124}\rput(4.3,1.65){\scriptsize 3124}
 \Cnode(2.5,3.3){1342}\rput(2.9,3.55){\scriptsize 1342}
 \Cnode(4.7,3.3){3241}\rput(4.35,3.55){\scriptsize 3241}
 \Cnode(2.5,4.7){4132}\rput(2.05,4.85){\scriptsize 4132}
 \Cnode(4.7,4.7){2431}\rput(5.15,4.85){\scriptsize 2431}
 \Cnode(3.6,4){1234}\rput(3.6,4.35){\scriptsize 1234}
 \Cnode(1.3,2.6){3412}\rput(0.8,2.6){\scriptsize 3412}
 \Cnode(5.9,2.6){4321}\rput(6.4,2.6){\scriptsize 4321}

\ncline{3412}{4213}\ncline{3412}{2314}\ncline{3412}{1342}\ncline{3412}{4132}
\ncline{4132}{2431}\ncline{4132}{1234}\ncline{4132}{1342}\ncline{1342}{3241}
\ncline{1342}{2143}\ncline{2314}{1234}\ncline{2314}{4213}\ncline{4213}{2143}
\ncline{4213}{1423}\ncline{1234}{2431}\ncline{1234}{3124}\ncline{2143}{1423}
\ncline{2143}{3241}\ncline{4321}{2431}\ncline{4321}{3241}\ncline{4321}{3124}
\ncline{4321}{1423}\ncline{3241}{2431}\ncline{3124}{1423}\ncline{2314}{3124}
\rput(3.6,0){\scriptsize $AG_4$}

\end{pspicture}
\caption{\label{f3}\footnotesize{Alternating group graphs $AG_3$ and
$AG_4$}}
\end{figure}
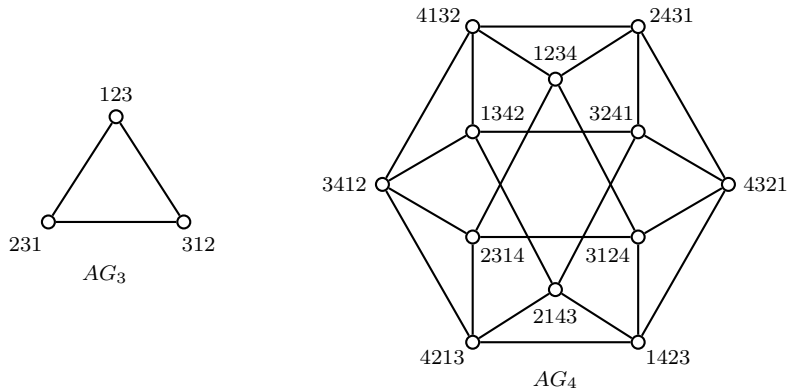

For $n\geqslant 3$, let $S=\{(1,2)(1,i), (1,i)(1,2):\ 3\leqslant
i\leqslant n\}$, where $(1,2)(1,i)$ and $(1,i)(1,2)$ are mutually
inverse for each $i=3,\cdots,n$, and so $S=S^{-1}$. The Cayley graph
$C_{\Gamma_n}(S)$ is called the {\em alternating group graph},
proposed by Jwo {\it et al.}~\cite{jhd93} in 1993 and denoted by
$AG_n$, which is $(2n-4)$-regular and $(2n-4)$-connected. $AG_{3}$
and $AG_{4}$ are shown in Fig.~\ref{f3}.

\vskip6pt

It is known that $\kappa_{2}(AG_n)=6n-19$ for $n\geqslant 5$
determined by Lin {\it et al.}~\cite{lzxw15} and $t_{c}(AG_4)=4$ and
$t_{c}(AG_n)=6n-19$ for $n\geqslant 6$ obtained by Zhou and
Xu~\cite{zx13}, and Hao {\it et al.}~\cite{hfz13}, in which
``$t_{c}(AG_n)=6n-18$" is a slip of the pen. Thus,
$t_{c}(AG_n)=6n-19=\kappa_{2}(AG_n)$ for $n\geqslant 6$.

\vskip10pt

\subsection{Star Graphs}

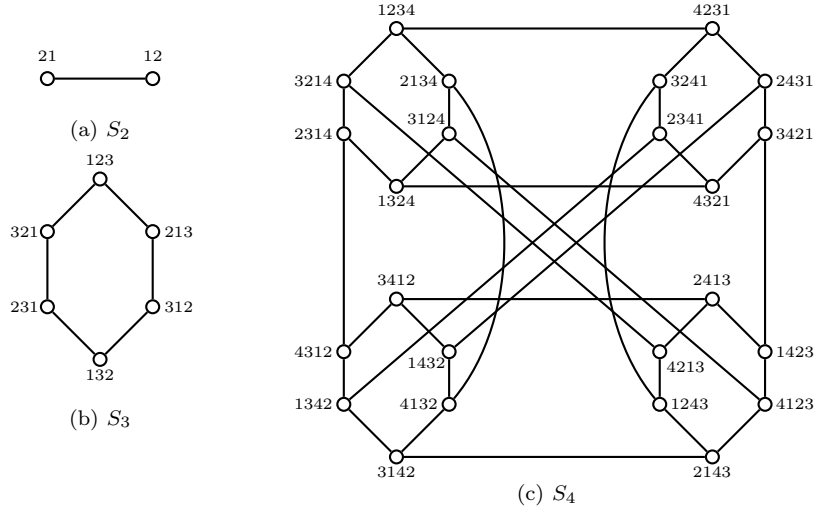
\begin{figure}[h]
\begin{pspicture}(-2,.3)(3,2.2)
\psset{radius=.1}

\Cnode(1.9,1){12}\rput(1.9,1.3){\tiny12}
\Cnode(.5,1){21}\rput(.5,1.3){\tiny21} \ncline{21}{12}
\rput(1.2,0.3){\scriptsize (a) $S_2$}
\end{pspicture}

\begin{pspicture}(-2,-2)(3,2.5)
\psset{radius=.1} \Cnode(.5,.5){231}\rput(.2,.5){\tiny231}
\Cnode(.5,1.5){321}\rput(.2,1.5){\tiny321}
\Cnode(1.9,1.5){213}\rput(2.25,1.5){\tiny213}
\Cnode(1.9,.5){312}\rput(2.25,.5){\tiny312}
\Cnode(1.2,-.2){132}\rput(1.2,-.4){\tiny132}
\Cnode(1.2,2.2){123}\rput(1.2,2.45){\tiny123}
\ncline{231}{321}\ncline{321}{123}\ncline{213}{123}
\ncline{213}{312}\ncline{312}{132}\ncline{132}{231}
\rput(1.2,-1){\scriptsize (b) $S_3$}
\end{pspicture}
\begin{pspicture}(-.5,0.5)(7.5,5.3)
\psset{radius=.1}

\Cnode(1.5,1){3142}\rput(1.5,.8){\tiny3142}
\Cnode(.8,1.7){1342}\rput(.4,1.7){\tiny1342}
\Cnode(2.2,1.7){4132}\rput(1.8,1.7){\tiny4132}
\Cnode(.8,2.4){4312}\rput(.4,2.4){\tiny4312}
\Cnode(2.2,2.4){1432}\rput(1.9,2.25){\tiny1432}
\Cnode(1.5,3.1){3412}\rput(1.5,3.35){\tiny3412}
\Cnode(1.5,4.6){1324}\rput(1.5,4.4){\tiny1324}
\Cnode(.8,5.3){2314}\rput(.4,5.3){\tiny2314}
\Cnode(2.2,5.3){3124}\rput(1.9,5.5){\tiny3124}
\Cnode(.8,6){3214}\rput(.4,6){\tiny3214}
\Cnode(2.2,6){2134}\rput(1.8,6){\tiny2134}
\Cnode(1.5,6.7){1234}\rput(1.5,6.95){\tiny1234}

\Cnode(5.7,1){2143}\rput(5.7,.8){\tiny2143}
\Cnode(5,1.7){1243}\rput(5.4,1.7){\tiny1243}
\Cnode(6.4,1.7){4123}\rput(6.8,1.7){\tiny4123}
\Cnode(5,2.4){4213}\rput(5.35,2.2){\tiny4213}
\Cnode(6.4,2.4){1423}\rput(6.8,2.4){\tiny1423}
\Cnode(5.7,3.1){2413}\rput(5.7,3.35){\tiny2413}
\Cnode(5.7,4.6){4321}\rput(5.7,4.4){\tiny4321}
\Cnode(5,5.3){2341}\rput(5.35,5.5){\tiny2341}
\Cnode(6.4,5.3){3421}\rput(6.8,5.3){\tiny3421}
\Cnode(5,6){3241}\rput(5.4,6){\tiny3241}
\Cnode(6.4,6){2431}\rput(6.8,6){\tiny2431}
\Cnode(5.7,6.7){4231}\rput(5.7,6.95){\tiny4231}

\ncline{3142}{4132}\ncline{3142}{1342}\ncline{4312}{1342}
\ncline{4312}{3412}\ncline{1432}{3412}\ncline{1432}{4132}
\ncline{1324}{2314}\ncline{1324}{3124}\ncline{3214}{2314}
\ncline{3214}{1234}\ncline{1234}{2134}\ncline{2134}{3124}\ncline{4312}{2314}
\ncline{2143}{4123}\ncline{2143}{1243}\ncline{1423}{4123}\ncline{2143}{3142}
\ncline{1423}{2413}\ncline{2413}{4213}\ncline{4213}{1243}
\ncline{4321}{2341}\ncline{4321}{3421}\ncline{3241}{2341}\ncline{3421}{1423}
\ncline{3241}{4231}\ncline{4231}{2431}\ncline{3421}{2431}\ncline{1234}{4231}
\ncline{3214}{4213}\ncline{3124}{4123}\ncline{2341}{1342}\ncline{1432}{2431}
\ncline{1324}{4321}\ncline{3412}{2413}
\nccurve[angleA=-50,angleB=50]{2134}{4132}
\nccurve[angleA=-130,angleB=130]{3241}{1243}
\rput(3.5,0.5){\scriptsize (c) $S_4$}

\end{pspicture}
\caption{\label{f4}\footnotesize{The star graphs $S_2, S_3$ and
$S_{4}$}}
\end{figure}

Let $\Omega_{n}$ be the symmetry group and $S=\{(1,i):\ 2\leqslant
i\leqslant n\}$. The Cayley graph $C_{\Omega_n}(S)$ is called a star
graph, denoted by $S_{n}$, proposed by Akers and
Krishnamurthy~\cite{ak89} in 1989. The graphs shown in
Figure~\ref{f4} are $S_2, S_3$ and $S_4$.

A star graph $S_{n}$ is $(n-1)$-regular and $(n-1)$-connected.
Furthermore, since a transposition changes the parity of a
permutation, each edge connects an odd permutation with an even
permutation, and so $S_{n}$ is bipartite, and contains no $C_4$. A
star graph is not only vertex-transitive but also
edge-transitive~\cite{ak89}, and so is symmetric.

\begin{lem}\label{lem3.1}
For any $x,y\in V(S_n)$, if $xy\notin E(S_n)$ and $N(x)\cap
N(y)\ne\emptyset$, then $|N(x)\cap N(y)|=1$.
\end{lem}

Since $S_n$ is $(n-1)$-regular and contains no $C_3$, according to
Lemma~\ref{lem3.1}, if $P_3=(x,y,z)$ is a 3-path, where $xz\notin
E(G)$, then $|N(x)\cap N(y)|=|N(y)\cap N(z)|=0$ and $N(x)\cap
N(z)=\{y\}$, and so the number of neighbors of $P_3$ in $S_n$ can be
counted as follows.
 $$
 \begin{array}{rl}
 |N(P_3)|&=d(x)+d(y)+d(z)-|N(x)\cap N(y)|-|N(y)\cap N(z)|-\Sigma(P_3)\\
  &=3(n-1)-4=3n-7.
  \end{array}
 $$
Since $S_n$ is vertex-transitive, for any $3$-path $P_3$ in $S_n$,
we have that
  \begin{equation}\label{e1}
  |N(P_3)|=3(n-1)-4=3n-7.
   \end{equation}

\begin{lem}\label{lem3.2} \textnormal{(Cheng and Lipt\'ak~\cite{cl06})}
Let $F\subset V(S_n)$ with $|F|\leqslant 3n-8$ and $n\geqslant 5$.
If $S_n-F$ is disconnected, then it has either two components, one
of which is an isolated vertex or an edge, or three components, two
of which are isolated vertices.
\end{lem}

Lin {\it et al.}~\cite{lthcl08}, Zhou and Xu~\cite{zx13} determined
$t_c(S_n)=3n-7$ for $n\geqslant 4$. However, $\kappa_{2}(S_n)$ has
not been determined so for. We can deduce these results by
Theorem~\ref{thm2.4}.

\begin{thm}\label{thm3.3}\
$t_{c}(S_n)=3n-7=\kappa_{2}(S_n)$ for $n\geqslant 5$.
\end{thm}

\begin{pf}
Since $S_n$ contains no $C_3$, $t=\min\{|N(T)|:\ T=P_3$ or $C_3$ in
$S_n\}=|N(P_3)|$, where $P_3$ is any 3-path in $S_n$ since $S_n$ is
vertex-transitive. Let $F=N(P_3)$. Then $|F|=t=3n-7$ by (\ref{e1}).
It is easy to check that $|V(S_n)|-|F|-3=n\,!-3n+4>0$ for
$n\geqslant 4$. Thus $F$ is a vertex-cut of $S_n$. To prove the
theorem, we only need to verify that $S_n$ satisfies conditions in
Theorem~\ref{thm2.4}.

{\rm (a)}\ If $|F|\leqslant t-1$ then, by Lemma~\ref{lem3.2},
$S_n-F$ has a large component and small components which contain at
most two vertices in total.

{\rm (b)}\ By Lemma~\ref{lem3.1}, $ \ell(S_n)=1$. Since $S_{n}$ is
$(n-1)$-regular bipartite, it contains no $5$-cycle, and so
$n-1\geqslant 4=2\ell(S_n)+2$.

{\rm (c)}\ When $n\geqslant 4$, it is easy to check that
 $$
 \begin{array}{rl}
 n\,!-n(t-1)-4&=n\,!-n(3n-8)-4\\
    &\geqslant 4(n-1)(n-2)-3n^2+8n-4\\
    &=(n-2)^2\\
    &>0.
 \end{array}
 $$

$S_n$ satisfies all conditions in Theorem~\ref{thm2.4}, and so
$t_c(S_n)=3n-7=\kappa_2(S_n)$.
\end{pf}

\vskip6pt

The star graph $S_n$ is an important topological structure of
interconnection networks and has attracted considerable attention
since it has been thought to be an attractive alternative to the
hypercube. However, since $S_n$ has $n\,!$ vertices, there is a
large gap between $n\,!$ and $(n+1)\,!$ for expanding $S_n$ to
$S_{n+1}$. To relax the restriction of the numbers of vertices in
$S_n$, the arrangement graph $A_{n,k}$ and the $(n,k)$-star graph
$S_{n,k}$ were proposed as generalizations of the star graph $S_n$.
In the following two sections, we discuss such two classes of
graphs, respectively.

For this purpose, we need some notations. Given two positive
integers $n$ and $k$ with $k<n$, let $P_{n,k}$ be a set of
arrangements of $k$ elements in $I_n$, i.e.,
$P_{n,k}=\{p_{1}p_{2}\ldots p_{k}:\ p_{i}\in I_n, p_{i}\neq p_{j},
1\leqslant i\neq j\leqslant k\}$. Clearly,
$|P_{n,k}|=\frac{n!}{(n-k)!}$.

\vskip10pt
%
\subsection{Arrangement Graphs}

The {\it $(n,k)$-arrangement graph}, denoted by $A_{n,k}$, was
proposed by Day and Tripathi~\cite{dt92} in 1992. The definition of
$A_{n,k}$ is as follows. $A_{n,k}$ has vertex-set $P_{n,k}$ and two
vertices are adjacent if and only if they differ in exactly one
position.

\vskip6pt

\begin{figure}[h]
\begin{pspicture}(-3.25,0)(7.5,6.5)
\psset{radius=.1}
\Cnode(.8,1){43}\rput(0.7,.7){\scriptsize 4\,3}
\Cnode(2.5,1){23}\rput(2.5,.7){\scriptsize 2\,3}
\Cnode(5.5,1){21}\rput(5.5,.7){\scriptsize 2\,1}
\Cnode(7.2,1){41}\rput(7.3,.7){\scriptsize 4\,1}
\Cnode(1.75,2.5){13}\rput(1.4,2.5){\scriptsize 1\,3}
\Cnode(4,2){24}\rput(4,1.7){\scriptsize 2\,4}
\Cnode(6.25,2.5){31}\rput(6.6,2.5){\scriptsize 3\,1}
\Cnode(3.3,3.5){14}\rput(3,3.7){\scriptsize 1\,4}
\Cnode(4.7,3.5){34}\rput(5.1,3.6){\scriptsize 3\,4}
\Cnode(4.7,4.7){32}\rput(5.1,4.7){\scriptsize 3\,2}
\Cnode(3.3,4.7){12}\rput(2.9,4.7){\scriptsize 1\,2}
\Cnode(4,6){42}\rput(4,6.3){\scriptsize 4\,2}

\ncline{43}{23}\ncline{21}{23}\ncline{41}{21}\ncline{43}{13}
\ncline{13}{23}\ncline{24}{23}\ncline{24}{21}\ncline{21}{31}
\ncline{13}{12}\ncline{13}{14}\ncline{42}{12}\ncline{14}{24}
\ncline{24}{34}\ncline{34}{31}\ncline{14}{34}\ncline{14}{12}
\ncline{32}{34}\ncline{32}{31}\ncline{12}{32}\ncline{12}{42}
\ncline{42}{32}\ncline{31}{41}
\nccurve[angleA=-160,angleB=100]{42}{43}
\nccurve[angleA=-20,angleB=80]{42}{41}
\nccurve[angleA=-30,angleB=-150]{43}{41} \rput(1.,5.5){ $A_{4,2}$}
\end{pspicture}
\caption{\label{f5}\footnotesize{The structure of a
$(4,2)$-arrangement graph $A_{4,2}$}}
\end{figure}
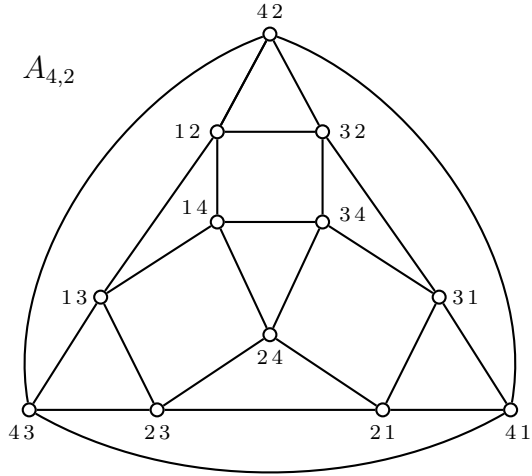

Figure~\ref{f5} shows a $(4,2)$-arrangement graph $A_{4,2}$, which
is isomorphic to $AG_4$ (see Fig.~\ref{f3}).

Since $|P_{n,k}|=\frac{n!}{(n-k)!}$ and $|S|=k(n-k)$, $A_{n,k}$ is a
$k(n-k)$-regular graph with order $\frac {n\,!}{(n-k)\,!}$, and is
$k(n-k)$-connected since $S$ is a minimal generating set of
$\Gamma_n$. Moreover, $A_{n,k}$ is vertex-transitive and
edge-transitive (see~\cite{dt92}), and so  $A_{n,k}$ is symmetric.
Clearly, $A_{n,1}\cong K_n$ and $A_{n,n-1}\cong S_n$. Chiang and
Chen~\cite{cc98} showed that $A_{n,n-2}\cong AG_n$. Thus, the
$(n,k)$-arrangement graph $A_{n,k}$ is naturally regarded as a
common generalization of the star graph $S_n$ and the alternating
group graph $AG_n$. For a fixed $i$ ($1\leqslant i\leqslant k)$, let
 $$
 V_i=\{p_1\cdots p_{i-1}q_ip_{i+1}\cdots p_k:\ q_i\in I_n\setminus\{p_1,
 \cdots, p_{i-1}, p_{i+1}, \cdots, p_k\}\}
$$
Then $|V_i|=n-k+1$. There are $|P_{n,k-1}|$ such $V_i$'s. By
definition, it is easy to see that the subgraph of $A_{n,k}$ induced
by $V_i$ is a complete graph $K_{n-k+1}$. In special,
$K_{n-k+1}=K_n$ if $k=1$, and $K_{n-k+1}=K_2$ if $k=n-1$.

When $n=k+1$, $A_{n,k}$ contains no $3$-cycle $C_3$, there is a big
difference in the way of dealing it with other conditions. Since
$A_{n,n-1}\cong S_n$, which has been discussed in the above
subsection, to avoid duplication of discussion, we may assume
$n\geqslant k+2$ and $k\geqslant 2$ in the following discussion.

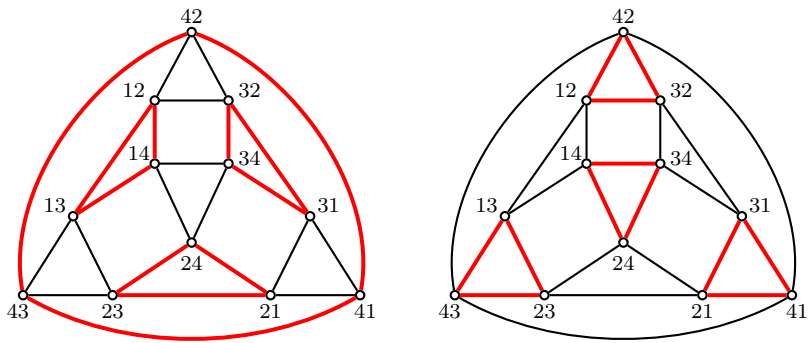
\begin{figure}[h]
\psset{unit=.7}
\begin{pspicture}(-2.5,0)(1.,6.5)
\psset{radius=.1} \Cnode(.8,1){43}\rput(0.7,.7){\scriptsize 43}
\Cnode(2.5,1){23}\rput(2.5,.7){\scriptsize 23}
\Cnode(5.5,1){21}\rput(5.5,.7){\scriptsize 21}
\Cnode(7.2,1){41}\rput(7.3,.7){\scriptsize 41}
\Cnode(1.75,2.5){13}\rput(1.4,2.7){\scriptsize 13}
\Cnode(4,2){24}\rput(4,1.6){\scriptsize 24}
\Cnode(6.25,2.5){31}\rput(6.6,2.7){\scriptsize 31}
\Cnode(3.3,3.5){14}\rput(3,3.7){\scriptsize 14}
\Cnode(4.7,3.5){34}\rput(5.1,3.6){\scriptsize 34}
\Cnode(4.7,4.7){32}\rput(5.1,4.9){\scriptsize 32}
\Cnode(3.3,4.7){12}\rput(2.9,4.9){\scriptsize 12}
\Cnode(4,6){42}\rput(4,6.3){\scriptsize 42}

 \ncline{43}{23}
 \ncline{41}{21}
 \ncline{43}{13}
 \ncline{13}{23}
 \ncline{21}{31}
 \ncline{14}{24}
 \ncline{24}{34}
 \ncline{14}{34}
 \ncline{12}{32}
 \ncline{12}{42}
 \ncline{42}{32}
 \ncline{31}{41}

 \ncline[linecolor=red,linewidth=1.5pt]{21}{23}
 \ncline[linecolor=red,linewidth=1.5pt]{24}{23}
 \ncline[linecolor=red,linewidth=1.5pt]{24}{21}

 \ncline[linecolor=red,linewidth=1.5pt]{13}{12}
 \ncline[linecolor=red,linewidth=1.5pt]{13}{14}
 \ncline[linecolor=red,linewidth=1.5pt]{34}{31}

 \ncline[linecolor=red,linewidth=1.5pt]{14}{12}
 \ncline[linecolor=red,linewidth=1.5pt]{32}{34}
 \ncline[linecolor=red,linewidth=1.5pt]{32}{31}

 \nccurve[linecolor=red,linewidth=1.5pt,angleA=-160,angleB=100]{42}{43}
 \nccurve[linecolor=red,linewidth=1.5pt,angleA=-20,angleB=80]{42}{41}
 \nccurve[linecolor=red,linewidth=1.5pt,angleA=-30,angleB=-150]{43}{41}
\end{pspicture}
\begin{pspicture}(-7.,0)(-3.5,6.5)
\psset{radius=.1} \Cnode(.8,1){43}\rput(0.7,.7){\scriptsize 43}
\Cnode(2.5,1){23}\rput(2.5,.7){\scriptsize 23}
\Cnode(5.5,1){21}\rput(5.5,.7){\scriptsize 21}
\Cnode(7.2,1){41}\rput(7.3,.7){\scriptsize 41}
\Cnode(1.75,2.5){13}\rput(1.4,2.7){\scriptsize 13}
\Cnode(4,2){24}\rput(4,1.6){\scriptsize 24}
\Cnode(6.25,2.5){31}\rput(6.6,2.7){\scriptsize 31}
\Cnode(3.3,3.5){14}\rput(3,3.7){\scriptsize 14}
\Cnode(4.7,3.5){34}\rput(5.1,3.6){\scriptsize 34}
\Cnode(4.7,4.7){32}\rput(5.1,4.9){\scriptsize 32}
\Cnode(3.3,4.7){12}\rput(2.9,4.9){\scriptsize 12}
\Cnode(4,6){42}\rput(4,6.3){\scriptsize 42}

 \ncline[linecolor=red,linewidth=1.5pt]{43}{23}
 \ncline{21}{23}
 \ncline[linecolor=red,linewidth=1.5pt]{41}{21}
 \ncline[linecolor=red,linewidth=1.5pt]{43}{13}
 \ncline[linecolor=red,linewidth=1.5pt]{13}{23}
 \ncline{24}{23}
 \ncline{24}{21}
 \ncline[linecolor=red,linewidth=1.5pt]{21}{31}
 \ncline{13}{12}
 \ncline{13}{14}
 \ncline{42}{12}
 \ncline[linecolor=red,linewidth=1.5pt]{14}{24}
 \ncline[linecolor=red,linewidth=1.5pt]{24}{34}
 \ncline{34}{31}
 \ncline[linecolor=red,linewidth=1.5pt]{14}{34}
 \ncline{14}{12}
 \ncline{32}{34}
 \ncline{32}{31}
 \ncline[linecolor=red,linewidth=1.5pt]{12}{32}
 \ncline[linecolor=red,linewidth=1.5pt]{12}{42}
 \ncline[linecolor=red,linewidth=1.5pt]{42}{32}
 \ncline[linecolor=red,linewidth=1.5pt]{31}{41}
 \nccurve[angleA=-160,angleB=100]{42}{43}
 \nccurve[angleA=-20,angleB=80]{42}{41}
 \nccurve[angleA=-30,angleB=-150]{43}{41}

\end{pspicture}

\caption{\label{f6}\footnotesize{Two partitions of $A_{4,2}$ into
$4$ triangles $K_3$ (red edges)}}
\end{figure}

Thus, when $n\geqslant k+2$ and $k\geqslant 2$, for each fixed $i$
($1\leqslant i\leqslant k)$, the vertex-set of $A_{n,k}$ can be
partitioned into $|P_{n,k-1}|$ subsets, each of which induces a
complete graph $K_{n-k+1}$. For example, for $n=4$ and $k=2$,
$|P_{4,1}|=4$. Fig.~\ref{f6} illustrates two partitions of
$V(A_{4,2})$ into $4$ subsets for each $i=1,2$, each of which
induces a complete graph $K_3$ (red edges). This fact and the
arbitrariness of $i\ (1\leqslant i\leqslant k)$ show that each
vertex is contained in $k$ distinct $K_{n-k+1}$'s, and each edge is
contained in $(n-k-1)$ distinct 3-cycles, that is, any two adjacent
vertices have exactly $(n-k-1)$ common neighbors.

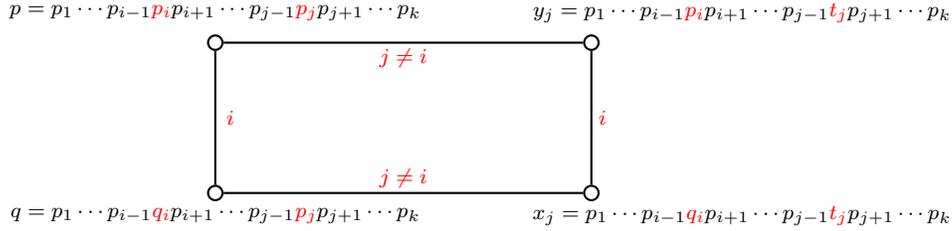
\begin{figure}[h]  
\psset{unit=1.}     
\begin{pspicture}(-1,-.3)(4,2.8)
 \cnode(3,2){.11}{1}\rput(3,2.4){\scriptsize$p=p_1\cdots p_{i-1}{\red p_i}p_{i+1}\cdots p_{j-1}{\red p_j}p_{j+1}\cdots p_k$}
 \cnode(8,2){.11}{2}\rput(10,2.4){\scriptsize$y_j=p_1\cdots p_{i-1}{\red p_i}p_{i+1}\cdots p_{j-1}{\red t_j}p_{j+1}\cdots p_k$}
 \cnode(8,0){.11}{3}\rput(3.,-.3){\scriptsize$q=p_1\cdots p_{i-1}{\red q_i}p_{i+1}\cdots p_{j-1}{\red p_j}p_{j+1}\cdots p_k$}
 \cnode(3,0){.11}{4}\rput(10,-0.3){\scriptsize$x_j=p_1\cdots p_{i-1}{\red q_i}p_{i+1}\cdots p_{j-1}{\red t_j}p_{j+1}\cdots p_k$}
 \ncline{1}{2}\ncline{1}{4} \ncline{2}{3} \ncline{3}{4}

 \rput(5.5,1.8){\scriptsize${\red j\ne i}$}\rput(5.5,0.2){\scriptsize${\red j\ne i}$}
 \rput(3.2,1){\scriptsize${\red i}$}\rput(8.15,1){\scriptsize${\red i}$}
\end{pspicture}

\caption{
\label{f7}                                       
\footnotesize  Construction of a chordless 4-cycle containing a
given edge $pq$ in $A_{n,k}$}
\end{figure}

Furthermore, each edge of $A_{n,k}$ is contained in $(k-1)$
chordless 4-cycles when $n\geqslant k+2$ and $k\geqslant 2$. In
fact, let $p\,q\in E(A_{n,k})$, if $p=p_1\cdots
p_{i-1}p_ip_{i+1}\cdots p_k$, then $q=p_1\cdots
p_{i-1}q_ip_{i+1}\cdots p_k$, where $q_i\in I_n\setminus\{p_1,
\cdots, p_k\}$. For each $j\in\{1,2,\cdots,k\}$ and $j\ne i$, let
 $$
 \begin{array}{rl}
 &x_j=p_1\cdots p_{i-1}q_ip_{i+1}\cdots p_{j-1}t_jp_{j+1}\cdots p_k\ \ {\rm and}\\
 &y_j=p_1\cdots p_{i-1}p_ip_{i+1}\cdots p_{j-1}t_jp_{j+1}\cdots p_k,
 \end{array}
 $$
where $t_j\in I_n\setminus\{p_1,\cdots, p_k,q_i\}$, such $t_j$
certainly exists since $n\geqslant k+2$ and $k\geqslant 2$. Then,
$(p,q,x_j,y_j)$ is a chordless 4-cycle in $A_{n,k}$ for each
$j\in\{1,2,\cdots,k\}$ and $j\ne i$ (see Fig.~\ref{f7}).

According to the above discussion, we have the following result.

\begin{lem}\label{lem3.4} When $n\geqslant k+2$,
for any $x,y\in V(A_{n,k})$, then $|N(x)\cap N(y)|=n-k-1$ if $xy\in
E(A_{n,k})$; $|N(x)\cap N(y)|\leqslant 2$ if $xy\notin E(A_{n,k})$
and $N(x)\cap N(y)\ne\emptyset$; and $|N(x)\cap N(y)|=0$ otherwise.
\end{lem}

Since each edge is contained in a $K_{n-k+1}$ ($n\geqslant k+2$),
for a 3-cycle $C_3=(x,y,z)$, every vertex in $V(K_{n-k+1}-C_3)$ is a
common neighbor of the three edges $xy, yz, zx$. In other words,
when we count the number $|N(C_3)|$ of neighbors of $C_3$ in
$A_{n,k}$, every vertex in $V(K_{n-k+1}-C_3)$ is counted three
times. Thus, the number $|N(C_3)|$ of neighbors of $C_3$ in
$A_{n,k}$ can be counted as follows.
 $$
 \begin{array}{rl}
 |N(C_3)|&=d(x)+d(y)+d(z)-2|V(K_{n-k+1}-C_3)|-\Sigma(C_3)\\
       &=3k(n-k)-2(n-k-2)-6\\
       &=(3k-2)(n-k)-2.
       \end{array}
 $$
Since $A_{n,k}$ is vertex-transitive, for any $3$-cycle $C_3$ in
$A_{n,k}$, we have that
 \begin{equation}\label{e2}
 |N(C_3)|=(3k-2)(n-k)-2.
 \end{equation}

Since $A_{n,k}$ contains chordless 4-cycle, say $(x,y,z,u)$, we
choose a $3$-path $P_3=(x,y,z)$. Then $xz\notin E(A_{n,k})$. Since
each edge is contained in a $K_{n-k+1}$, $|N(x)\cap N(y)|=|N(y)\cap
N(z)|=n-k-1$ and $|N(z)\cap N(x)|=|\{y,u\}|=2$ by
Lemma~\ref{lem3.4}. Note that two edge $xy$ and $yz$ are in
different complete graphs. Thus, the number of neighbors of $P_3$ in
$A_{n,k}$ can be counted as follows.
 $$
 \begin{array}{rl}
 |N(P_3)|&=d(x)+d(y)+d(z)-|N(x)\cap N(y)|\\
         &\quad -|N(y)\cap N(z)|-|N(z)\cap N(x)\setminus\{y\}|-\Sigma(P_3)\\
       &=3k(n-k)-2(n-k-1)-1-4\\
       &=(3k-2)(n-k)-3
       \end{array}
 $$
Since $A_{n,k}$ is vertex-transitive, for any $3$-path $P_3$ in
$A_{n,k}$, we have that
 \begin{equation}\label{e3}
  |N(P_3)|=(3k-2)(n-k)-3
 \end{equation}

\begin{lem}{\rm\cite{zx13}}\label{lem3.5}
Let $F$ be a vertex-cut of $A_{n,k}$ with $|F|\leqslant
(3k-2)(n-k)-4$. If $n\geqslant k+2$ and $k\geqslant 4$, then
$A_{n,k}-F$ contains either two components, one of which is an
isolated vertex or an isolated edge, or three components, two of
which are isolated vertices.
\end{lem}

Zhou and Xu~\cite{zx13} determined that for $n\geqslant k+2$ and
$k\geqslant 4$, $t_{c}(A_{n,k})=(3k-2)(n-k)-3$. However,
$\kappa_2(A_{n,k})$ has not been determined. We can deduce these
results by Theorem~\ref{thm2.4}.

\begin{thm}\label{thm3.6}\
$t_{c}(A_{n,k})=(3k-2)(n-k)-3=\kappa_{2}(A_{n,k})$ for $n\geqslant
k+2$ and $k(n-k)\geqslant 8$.
\end{thm}

\begin{pf}
Comparing (\ref{e2}) with (\ref{e3}), when $n\geqslant k+2$,
$t=\min\{|N(T)|:\ T=P_3$ or $C_3$ in $A_{n,k}\}=|N(P_3)|$, where
$P_3$ is any 3-path in $A_{n,k}$ since $A_{n,k}$ is
vertex-transitive. Let $F=N(P_3)$. Then $|F|=t=(3k-2)(n-k)-3$ by
(\ref{e3}). It is easy to check that $F$ is a vertex-cut of
$A_{n,k}$. To prove the theorem, we only need to verify that
$A_{n,k}$ satisfies conditions in Theorem~\ref{thm2.4}.

{\rm (a)}\ If $|F|\leqslant t-1$ then, by Lemma~\ref{lem3.5},
$A_{n,k}-F$ has a large component and small components which contain
at most two vertices in total.

{\rm (b)}\ By Lemma~\ref{lem3.4}, $\ell(A_{n,k})=2$, and so
 $k(n-k)\geqslant 8=3\ell(A_{n,k})+2$.

 {\rm (c)}\ It is not difficult to check that
 $$
 \begin{array}{rl}
 &\quad |V|-[(\Delta+1)(t-1)+4]\\
 &= |V|-(k(n-k)+1)((3k-2)(n-k)-4)-4\\
    &=|V|-3k^2(n-k)^2+2k(n-k)^2+(k+2)(n-k)\\
    &>|V|-3k^2(n-k)^2\quad \text{(for $n-k\geqslant 2$)}\\
    &\geqslant |V|-3(n-2)^2(n-k+1)^2\quad \text{(for $k\leqslant n-2$)}\\
    &=n!/(n-k)!-3(n-2)^2(n-k+1)^2\\
    &=n(n-1)\cdots(n-k+1)-3(n-2)^2(n-k+1)^2\\
    &>3(n-2)^2(n-k+1)^2-3(n-2)^2(n-k+1)^2\\
    &=0.
 \end{array}
 $$
$A_{n,k}$ satisfies all conditions in Theorem~\ref{thm2.4}, and so
$t_c(A_{n,k})=(3k-2)(n-k)-3=\kappa_2(A_{n,k})$.
\end{pf}


Since $A_{n,n-2}\cong AG_n$, by Theorem~\ref{thm3.6}, we immediately
obtain the following results.

\begin{cor}\
$t_{c}(AG_n)=6n-19=\kappa_{2}(AG_n)$ for $n\geqslant 6$.
\end{cor}

\vskip10pt

\subsection{$(n,k)$-Star Graphs}

The $(n,k)$-star graph $S_{n,k}$, proposed by Chiang {\it et
al.}~\cite{cc95} in 1995 as another generalization of the star graph
$S_n$, has vertex-set $P_{n,k}$, a vertex $p=p_{1}p_{2}\ldots
p_{i}\ldots p_{k}$ is adjacent to a vertex

{\rm (a)}\ $p_{i}p_{2}\cdots p_{i-1}p_{1}p_{i+1}\cdots p_{k}$, where
$i\in\{2,3,\cdots,k\}$ (swap-edge).

{\rm (b)} $p'_{1} p_{2}p_{3}\cdots p_{k}$, where $p'_{1}\in
I_n\setminus \{p_{i}:\ i \in I_k \}$ (unswap-edge).

\vskip6pt

\begin{figure}[h]  
\psset{unit=1.}     
\begin{pspicture}(-3.5,0.8)(0,.7)
\rput{0}{%
\SpecialCoor\degrees[6]
\multido{\i=0+1}{6}{\rput(1;\i){\cnode{.11}{1\i}}}
\ncline{10}{11}\ncline{11}{12}\ncline{12}{13}\ncline{13}{14}\ncline{14}{15}\ncline{15}{10}
}%
\end{pspicture}
\begin{pspicture}(-4.,0.65)(0,.7)
\rput{0}{%
\SpecialCoor\degrees[6]
\multido{\i=0+1}{6}{\rput(1;\i){\cnode{.11}{2\i}}}
\ncline{20}{21}\ncline{21}{22}\ncline{22}{23}\ncline{23}{24}\ncline{24}{25}\ncline{25}{20}
}%
\end{pspicture}

\begin{pspicture}(-3.5,-1.)(0,3.)
\rput{0}{%
\SpecialCoor\degrees[6]
\multido{\i=0+1}{6}{\rput(1;\i){\cnode{.11}{3\i}}}
\ncline{30}{31}\ncline{31}{32}\ncline{32}{33}\ncline{33}{34}\ncline{34}{35}\ncline{35}{30}
}%
\end{pspicture}
\begin{pspicture}(-4.,-1.)(0,3.)
\rput{0}{%
\SpecialCoor\degrees[6]
\multido{\i=0+1}{6}{\rput(1;\i){\cnode{.11}{4\i}}}
\ncline{40}{41}\ncline{41}{42}\ncline{42}{43}\ncline{43}{44}\ncline{44}{45}\ncline{45}{40}
}%

\ncline{11}{22}\rput(-3.6,5.){\scriptsize$431$}\rput(-.6,5.){\scriptsize$134$}
\ncline{35}{44}\rput(-3.6,-1.15){\scriptsize$243$}\rput(-.6,-1.15){\scriptsize$342$}
\ncline{13}{33}\rput(-4.7,3.8){\scriptsize$321$}\rput(-4.7,-0.){\scriptsize$123$}
\ncline{10}{30}\rput(-2.7,3.8){\scriptsize$341$}\rput(-2.7,-0.){\scriptsize$143$}
\ncline{20}{40}\rput(-1.4,3.8){\scriptsize$234$}\rput(-1.4,-0.){\scriptsize$432$}
\ncline{23}{43}\rput(0.6,3.8){\scriptsize$214$}\rput(0.6,-0.){\scriptsize$412$}
\ncline{12}{42}\rput(-4.6,5.){\scriptsize$231$}\rput(-.3,0.55){\scriptsize$132$}
\ncline{15}{45}\rput(-3.8,3.3){\scriptsize$241$}\rput(.6,-1.15){\scriptsize$142$}
\ncline{31}{21}\rput(-3.85,0.6){\scriptsize$413$}\rput(.6,5.){\scriptsize$314$}
\ncline{34}{24}\rput(-4.6,-1.15){\scriptsize$423$}\rput(-.3,3.25){\scriptsize$324$}

\rput(-4.7,1.15){\scriptsize$213$}\rput(.55,1.15){\scriptsize$312$}

\rput(-4.7,2.7){\scriptsize$213$}\rput(.6,2.7){\scriptsize$421$}

\psarc(-2.1,-5.25){2.62in}{68}{112}
\psarc(-2.1,9.1){2.62in}{248}{292}
 \rput(-1.9,-1.4){\scriptsize $S_{4,3}$}
\end{pspicture}
\begin{pspicture}(-.4,-.2)(3,5.3)
\psset{radius=.12}

\Cnode(2.2,0.5){32}\rput(2.,.2){\scriptsize32}
\Cnode(5.,0.5){23}\rput(5.2,0.2){\scriptsize23}
\Cnode(1.4,1.6){12}\rput(1.7,1.8){\scriptsize12}

\Cnode(5.8,1.6){13}\rput(5.5,1.8){\scriptsize 13}
\Cnode(2.7,1.6){42}\rput(3.,1.5){\scriptsize42}
\Cnode(4.6,1.6){43}\rput(4.3,1.5){\scriptsize43}

\Cnode(2.8,4.9){21}\rput(2.5,4.9){\scriptsize 21}
\Cnode(4.3,4.9){31}\rput(4.65,4.9){\scriptsize 31}
\Cnode(3.6,4){41}\rput(3.9, 3.9){\scriptsize 41}

\Cnode(3.6,3.2){14}\rput(3.9,3.2){\scriptsize 14}
\Cnode(3.,2.3){24}\rput(2.7,2.4){\scriptsize24}
\Cnode(4.2,2.3){34}\rput(4.5,2.4){\scriptsize 34}

\ncline{32}{42}\ncline{42}{12}\ncline{12}{32}
\ncline{23}{13}\ncline{13}{43}\ncline{43}{23}
\ncline{21}{31}\ncline{31}{41}\ncline{41}{21}
\ncline{14}{24}\ncline{24}{34}\ncline{34}{14}

\ncline{32}{23}\ncline{13}{31}\ncline{21}{12}
\ncline{14}{41}\ncline{24}{42}\ncline{43}{34}

\rput(3.6,-.2){\scriptsize $S_{4,2}$}

\end{pspicture}

\vskip.4cm \caption{\label{f8}\footnotesize {Two $(n,k)$-star graphs
$S_{4,3}$ and $S_{4,2}$}}
\end{figure}
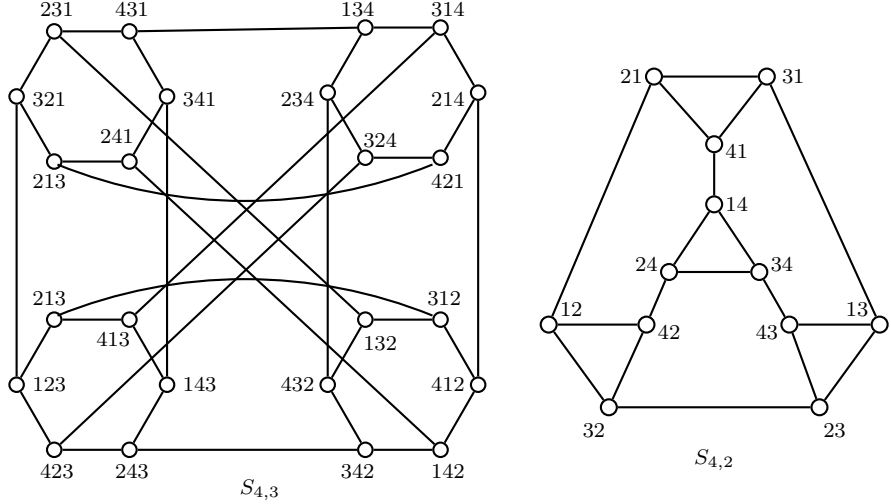

Figure~\ref{f8} shows two $(n,k)$-star graphs $S_{4,3}$ and
$S_{4,2}$, where $S_{4,3}\cong S_4$ and $S_{4,2}\cong AN_4$.

\vskip16pt

Since $|P_{n,k}|=\frac{n!}{(n-k)!}$ and $|S|=n-1$, $S_{n,k}$ is an
$(n-1)$-regular and $(n-1)$-connected graph with order $\frac
{n\,!}{(n-k)\,!}$. Moreover, $S_{n,k}$ is vertex-transitive,
however, it is not edge-transitive if $n\geqslant k+2$ (see Chiang
{\it et al.}~\cite{cc95}).

By definition, $S_{n,1}\cong K_n$ and $S_{n,n-1}\cong S_n$
obviously. Moreover, Cheng {\it et al.}~\cite{cqs12} showed
$S_{n,n-2}\cong AN_n$. Thus, the $(n,k)$-star graph $S_{n,k}$ is
naturally regarded as a common generalization of the star graph
$S_n$ and the alternating group network $AN_n$.
%

For any $\alpha=p_2p_3\cdots p_{k}\in P_{n,k-1}$ $(2\leqslant k
\leqslant n)$, let
 $$
 V_\alpha=\{p_1\alpha:\ p_1\in I_n\setminus \{p_i:\ 2\leqslant i\leqslant k\}.
 $$
By definition, it is easy to see that the subgraph of $S_{n,k}$
induced by $V_\alpha$ is a complete graph $K_{n-k+1}$. Thus,
$V(S_{n,k})$ can be partitioned into $|P_{n,k-1}|$ subsets, each of
which induces a complete graph $K_{n-k+1}$ whose edges are
unswap-edges. Furthermore, there is at most one swap-edge between
any two complete graphs, and so $S_{n,k}$ contains neither $4$-cycle
nor 5-cycle.

\begin{lem}\label{lem3.8} {\rm \cite{lx14}}
For any $x,y\in V(S_{n,k})$, then $|N(x)\cap N(y)|=n-k-1$ if $xy\in
E(S_{n,k})$ is an unswap-edge, $|N(x)\cap N(y)|=1$ if $xy\notin
E(S_{n,k})$ and $N(x)\cap N(y)\ne\emptyset$, and $|N(x)\cap N(y)|=0$
otherwise.
\end{lem}

Since $K_{n-k+1}=K_n$ when $k=1$ and $K_{n-k+1}=K_2$ when $k=n-1$,
like $A_{n,k}$, to avoid duplication of discussion, we may assume
$n\geqslant k+2$ and $k\geqslant 2$ in the following discussion.

For a 3-cycle $C_3=(x,y,z)$, since it is contained in a complete
graph $K_{n-k+1}$, every vertex in $V(K_{n-k+1}-C_3)$ is a common
neighbor of the tree edges $xy, yz, zx$. In other words, when we
count the number of neighbors of $C_3$ in $S_{n,k}$, every vertex in
$V(K_{n-k+1}-C_3)$ is counted three times. Thus, the number of
neighbors of $C_3$ in $S_{n,k}$ can be counted as follows.
 $$
 \begin{array}{rl}
 |N(C_3)|&=d(x)+d(y)+d(z)-2|V(K_{n-k+1}-C_3)|-\Sigma(C_3)\\
       &=3(n-1)-2(n-k-2)-6\\
       &=n+2k-5.
       \end{array}
 $$
Since $S_{n,k}$ is vertex-transitive, for any $3$-cycle $C_3$ in
$S_{n,k}$, we have that
 \begin{equation}\label{e4}
 |N(C_3)|=n+2k-5.
  \end{equation}

For a $3$-path $P_3=(x,y,z)$ with $xz\notin E(S_{n,k})$, then one of
two edges $xy$ and $yz$ is an unswap-edge and another is a
swap-edge. Without loss of generality, suppose that $xy$ is an
unswap-edge and $yz$ is a swap-edge. Then $|N(x)\cap N(y)|=n-k-1$,
$|N(y)\cap N(z)|=0$ and $|N(z)\cap N(x)|=|\{y\}|=1$ by
Lemma~\ref{lem3.8}. Thus, the number of neighbors of $C_3$ in
$A_{n,k}$ can be counted as follows.
 $$
 \begin{array}{rl}
 |N(P_3)|&=d(x)+d(y)+d(z)-|N(x)\cap N(y)|\\
         &\quad -|N(y)\cap N(z)|-|N(z)\cap N(x)\setminus\{y\}|-\Sigma(P_3)\\
       &=3(n-1)-(n-k-1)-0-4\\
       &=2n+k-6.
       \end{array}
 $$
Since $S_{n,k}$ is vertex-transitive, for any $3$-path $P_3$ in
$S_{n,k}$, we have that
 \begin{equation}\label{e5}
  |N(P_3)|=2n+k-6.
   \end{equation}

\begin{lem}\label{lem3.9} {\rm \cite{z12}}
Let $F$ be a vertex-cut of $S_{n,k}$ ($n\geqslant k+2$ and
$k\geqslant 3$) with $|F|\leqslant n+2k-6$. Then $S_{n,k}-F$
contains either two components, one of which is an isolated vertex
or an isolated edge, or three components, two of which are both
isolated vertices.
\end{lem}

Zhou~\cite{z12} determined that $t_{c}(S_{n,k})= n+2k-5$ if
$n\geqslant k+2$ and $k\geqslant 3$. However, $\kappa_2(S_{n,k})$
has not been determined. We can deduce these results by
Theorem~\ref{thm2.4}.

\begin{thm}\label{thm3.10}\
$t_{c}(S_{n,k})= n+2k-5=\kappa_{2}(S_{n,k})$ if $n\geqslant k+2$ and
$k\geqslant 3$.
\end{thm}

\begin{pf}
Let $t=\min\{|N(T)|:\ T=P_3$ or $C_3$ in $S_{n,k}\}$. By
Lemma~\ref{lem3.8}, $S_{n,k}$ contains 3-cycles when $n\geqslant
k+2$. Comparing (\ref{e4}) with (\ref{e5}), $t=|N(C_3)|=n+2k-5$,
where $C_3$ is any 3-cycle in $S_{n,k}$ since $S_{n,k}$ is
vertex-transitive. Let $F=N(C_3)$. Then $|F|=t$ and $F$ is a
vertex-cut of $S_{n,k}$. To prove the theorem, we only need to
verify that $S_{n,k}$ satisfies conditions in Theorem~\ref{thm2.4}.

{\rm (a)}\ If $|F|\leqslant t-1$ then, by Lemma~\ref{lem3.9},
$S_{n,k}-F$ has a large component and small components which contain
at most two vertices in total.

{\rm (b)}\ Since $S_{n,k}$ is $(n-1)$-regular and contains no
$5$-cycle $C_5$, by Lemma~\ref{lem3.8}, $\ell(S_{n,k})=1$, and so
$n-1\geqslant 4=2\ell(S_{n,k})+2$.

 {\rm (c)}\ It is not difficult to check that
 $$
 \begin{array}{rl}
 |V|-[n(t-1)+4]&=|V|-n(n+2k-6)-4\\
 &\geqslant |V|-n(3n-10)-4\quad \text{(for $k\leqslant n-2$)}\\
 &\geqslant |V|-3n(n-3)\quad \text{(for $n\geqslant 5$)}\\
 &\geqslant n(n-1)(n-2)-3n(n-3)\\
 &>3n(n-3)-3n(n-3)\\
 &=0.
 \end{array}
 $$

$S_{n,k}$ satisfies all conditions in Theorem~\ref{thm2.4}, and so
$t_c(S_{n,k})= n+2k-5=\kappa_2(S_{n,k})$. The theorem follows.
\end{pf}

\vskip6pt

Since $S_{n,n-2}\cong AN_n$, by Theorem~\ref{thm3.10}, we
immediately obtain the following results.

\begin{cor}
$t_{c}(AN_n)=3n-9=\kappa_{2}(AN_n)$ for $n\geqslant 5$.
\end{cor}

\vskip10pt

\subsection{Transposition Graphs}

Let $\mathscr T_n$ be a set of transpositions from $\Omega_n$ and
$S\subseteq\mathscr T_n$. The graph $T_S$ with vertex-set $I_n$ and
edge-set $\{ij:\ (i, j)\in S\}$ is called the {\it transposition
generating graph} or simply {\it transposition graph}. The Cayley
graph $C_{\Omega_n}(S)$ on $\Omega_n$ with respect to $S$ has $n\,!$
vertices.

For example, if $S=\{(1,i):\ 2\leqslant i\leqslant n\}$, then $T_S$
is a star $K_{1,n-1}$, the corresponding Cayley graph
$C_{\Omega_n}(S)$ is a star graph $S_n$, proposed by Akers and
Krishnamurthy~\cite{ak89}, perhaps, this is why they called such a
graph for the star graph.

\begin{figure}[h]
\begin{center}
\psset{unit=1pt}
\begin{pspicture}(0,10)(260,180)

\cnode(0,140){3}{a} \cnode(30,140){3}{b} \ncline{a}{b}

\rput(-12,140){{\tiny $12$}} \rput(42,140){{\tiny $21$}}

\rput(15,130){{\scriptsize $B_2$}}


\cnode(15,30){3}{a1} \cnode(15,78){3}{a4} \cnode(-6,42){3}{a6}
\cnode(-6,66){3}{a5} \cnode(36,42){3}{a2} \cnode(36,66){3}{a3}
\ncline{a1}{a2} \ncline{a2}{a3} \ncline{a3}{a4} \ncline{a4}{a5}
\ncline{a5}{a6} \ncline{a6}{a1}

\rput(48,40){{\tiny $231$}} \rput(48,68){{\tiny $213$}}
\rput(-17,40){{\tiny $312$}} \rput(-17,68){{\tiny $132$}}
\rput(15,87){{\tiny $123$}} \rput(15,21){{\tiny $321$}}

\rput(15,5){{\scriptsize $B_3$}}



\cnode(120,74){3}{b6} \cnode(120,98){3}{b5} \cnode(141,110){3}{b4}
\cnode(141,62){3}{b1} \cnode(162,98){3}{b3} \cnode(162,74){3}{b2}
\ncline{b1}{b2} \ncline{b2}{b3} \ncline{b3}{b4} \ncline{b4}{b5}
\ncline{b5}{b6} \ncline{b6}{b1}

\cnode(186,74){3}{c6} \cnode(186,98){3}{c5} \cnode(207,110){3}{c4}
\cnode(207,62){3}{c1} \cnode(228,98){3}{c3} \cnode(228,74){3}{c2}
\ncline{c1}{c2} \ncline{c2}{c3} \ncline{c3}{c4} \ncline{c4}{c5}
\ncline{c5}{c6} \ncline{c6}{c1}

\cnode(219,17){3}{d6} \cnode(219,41){3}{d5} \cnode(240,53){3}{d4}
\cnode(240,5){3}{d1} \cnode(261,41){3}{d3} \cnode(261,17){3}{d2}
\ncline{d1}{d2} \ncline{d2}{d3} \ncline{d3}{d4} \ncline{d4}{d5}
\ncline{d5}{d6} \ncline{d6}{d1}

\cnode(219,131){3}{e6} \cnode(219,155){3}{e5} \cnode(240,167){3}{e4}
\cnode(240,119){3}{e1} \cnode(261,155){3}{e3} \cnode(261,131){3}{e2}
\ncline{e1}{e2} \ncline{e2}{e3} \ncline{e3}{e4} \ncline{e4}{e5}
\ncline{e5}{e6} \ncline{e6}{e1}

\ncline{b2}{c6} \ncline{b3}{c5} \ncline{d5}{c1} \ncline{d4}{c2}
\ncline{e1}{c3} \ncline{e6}{c4}

\ncline{b1}{d6} \ncline{e2}{d3} \ncline{e5}{b4}


\psarc(207,86){88}{-50}{50} \psarc(207,86){88}{70}{170}
\psarc(207,86){88}{190}{290}


\rput(140,5){{\scriptsize $B_4$}}

\rput(107,74){{\tiny $4312$}} \rput(107,98){{\tiny $4132$}}
\rput(141,119){{\tiny $1432$}} \rput(141,53){{\tiny $3412$}}
\rput(162,108){{\tiny $1342$}} \rput(162,64){{\tiny $3142$}}

\rput(186,64){{\tiny $3124$}} \rput(186,108){{\tiny $1324$}}
\rput(200,120){{\tiny $1234$}} \rput(200,52){{\tiny $3214$}}
\rput(242,98){{\tiny $2134$}} \rput(242,74){{\tiny $2314$}}

\rput(214,9){{\tiny $3421$}} \rput(208,40){{\tiny $3241$}}
\rput(252,56){{\tiny $2341$}} \rput(250,0){{\tiny $4321$}}
\rput(273,44){{\tiny $2431$}} \rput(275,17){{\tiny $4231$}}

\rput(207,133){{\tiny $1243$}} \rput(215,163){{\tiny $1423$}}
\rput(240,177){{\tiny $4123$}} \rput(252,117){{\tiny $2143$}}
\rput(273,159){{\tiny $4213$}} \rput(273,128){{\tiny $2413$}}

\end{pspicture}

\end{center}
\caption{\label{f9}\footnotesize{The bubble-sort graphs $B_2$, $B_3$
and $B_4$. }}
\end{figure}
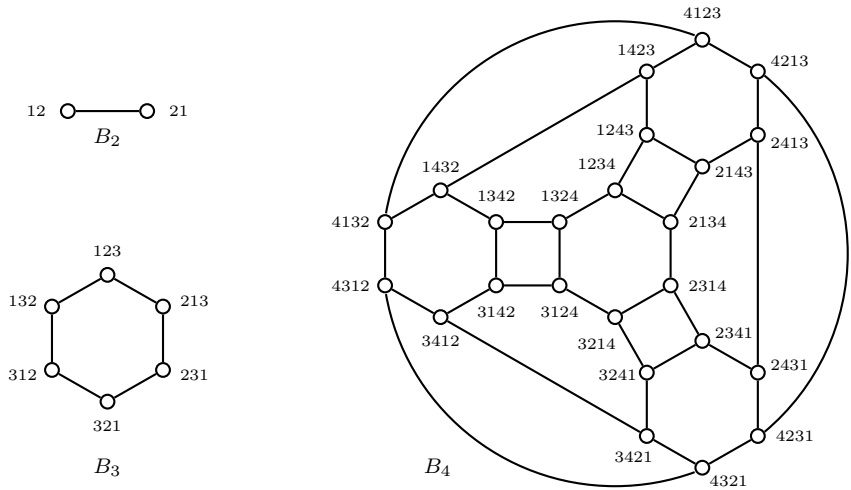

Here is another example, if $S=\{(i,i+1):\ 1\leqslant i\leqslant
n-1\}$, then $T_S$ is an $n$-path $P_n$, the corresponding Cayley
graph $C_{\Omega_n}(S)$ is called a bubble-sort graph $B_n$,
proposed by Akers and Krisnamurthy~\cite{ak89} in 1989. This series
of transpositions looks like to be along a straight line on the
bubbled. Perhaps this is why Akers and Krisnamurthy called such a
graph for the bubble-sort graph. Figure~\ref{f9} shows the
bubble-sort graphs $B_2$, $B_3$ and $B_4$.

It is a well-known result, due to Polya (see Berge~\cite{b71},
p118)), that a set $S\subseteq\mathscr T_n$ with $|S|=(n-1)$
generates $\Omega_n$ if and only if the transposition graph $T_S$ is
a tree, called a {\it transposition tree}.

Thus, one is interested in such a Cayley graph $C_{\Omega_n}(S)$
obtained from a transposition generating tree $T_S$, denoted by
$\mathscr T_n(S)$ shortly. The Cayley graph $\mathscr T_n(S)$ is a
bipartite graph since a transposition changes the parity of a
permutation, each edge connects an odd permutation with an even
permutation.

As we have seen from the above examples, $\mathscr T_n(S)$ is a star
graph $S_n$ if $T_S\cong K_{1,n-1}$, and a bubble-sort graph $B_n$
if $T_S\cong P_n$. Thus, the star graph $S_n$ and the bubble-sort
graph $B_n$ are special cases of the Cayley graph $\mathscr T_n(S)$.

Since when $T_S\cong K_{1,n-1}$, $\mathscr T_n(S)$ is a star graph
$S_n$. To avoid duplication of discussion, we may assume that $T_S$
is not a star $K_{1,n-1}$ in the following discussion.

Under this assumption, when $n\geqslant 4$, Lin {\it et
al.}~\cite{lthcl08} determined $t_c(\mathscr T_n(S))=3n-8$, Yang
{\it et al.}~\cite{ylm10} determined $\kappa_2(\mathscr
T_n(S))=3n-8$. We can deduce these results for $n\geqslant 7$ by
Theorem~\ref{thm2.4}.

According to the recursive architecture of $\mathscr T_n(S)$, we
easy obtain the following lemma.

\begin{lem}\label{lem3.12}
For any $x,y\in V(\mathscr T_n(S))$, if $xy\notin E(\mathscr
T_n(S))$ and $N(x)\cap N(y)\ne\emptyset$, then $|N(x)\cap N(y)|=1$
if $\mathscr T_n(S)=S_n$, and $|N(x)\cap N(y)|\leqslant 2$
otherwise.
\end{lem}

\begin{lem}\label{lem3.13} {\rm (Cheng and Lipt\'ak~\cite{cl06})}
For $n\geqslant 5$, if $T\subset V(\mathscr T_n(S))$ is a vertex-cut
with $|T|\leqslant 3n-8$, then $\mathscr T_n(S)-T$ contains either
two components, one of which is an isolated vertex or an isolated
edge, or three components, two of which are both isolated vertices.
%
%
%
%
%
%
\end{lem}

\begin{thm}\label{thm3.14}\
$t_{c}(\mathscr T_n(S))=3n-8=\kappa_{2}(\mathscr T_n(S))$ for
$n\geqslant 7$.
\end{thm}

\begin{pf}
Since $\mathscr T_n(S)$ is a partite graph, it contains no $C_3$,
and so $t=\min\{|N(T)|:\ T$ is a $3$-path or a $3$-cycle in
$\mathscr T_n(S)\}=|N(P_3)|$, where $P_3$ is any 3-path in $\mathscr
T_n(S)$ since $\mathscr T_n(S)$ is vertex-transitive. When $\mathscr
T_n(S)$ is not a star graph, it contains $C_4$, and so
$t=|N(P_3)|=3(n-1)-1-4=3n-8$. Let $F=N(P_3)$. It is easy to check
that $F$ is a vertex-cut of $\mathscr T_n(S)$. To prove the theorem,
we only need to verify that $\mathscr T_n(S)$ satisfies conditions
in Theorem~\ref{thm2.4}.

{\rm (a)}\ If $|F|\leqslant t-1$ then, by Lemma~\ref{lem3.13},
$\mathscr T_n(S)-F$ has a large component and small components have
at most two vertices in total.

{\rm (b)}\ By Lemma~\ref{lem3.12}, if $\mathscr T_n(S)\ne S_n$, then
$\ell(\mathscr T_n(S))=2$. Since $\mathscr T_n(S)$ is a bipartite
graph, it contains no $5$-cycle $C_5$. It follows that $
n-1\geqslant 6=2\ell(A_{n,k})+2$.

 {\rm (c)}\ It is easy to check that
$n!-[n(t-1)+4]>0$.

$\mathscr T_n(S)$ satisfies all conditions in Theorem~\ref{thm2.4},
and so $t_c(\mathscr T_n(S))=3n-8=\kappa_2(\mathscr T_n(S))$.
\end{pf}

\vskip6pt

Since when $T_S\cong P_n$ the Cayley graph $C_{\Omega_n}(S)$ is a
bubble-sort graph $B_n$, by Theorem~\ref{thm3.14}, we immediately
obtain the following result.

\begin{cor}\
$t_{c}(B_n)=3n-8=\kappa_{2}(B_n)$ for $n\geqslant 7$.
\end{cor}

 \vskip20pt

\subsection{$k$-ary $n$-cube Networks}

We first introduce the Cartesian product of graphs.

Let $G_1=(V_1,E_1)$ and $G_2=(V_2,E_2)$ be two undirected graphs.
The {\it Cartesian product} of $G_1$ and $G_2$ is an undirected
graph, denoted by $G_1\times  G_2$, where $V(G_1\times
G_2)=V_1\times  V_2$, two distinct vertices $x_1x_2$ and $y_1y_2$,
where $x_1,y_1\in V(G_1)$ and $x_2, y_2\in V(G_2)$, are linked by an
edge in $G_1\times G_2$ if and only if either $x_1=y_1$ and
$x_2y_2\in E(G_2)$, or $x_2=y_2$ and $x_1y_1\in E(G_1)$.

\begin{figure}[h]
\begin{pspicture}(-2.2,-.5)(1,3)
\cnode(1,1){.1}{0}\rput(.75,1){\scriptsize0}
\cnode(1,3){.1}{1}\rput(.75,3){\scriptsize1}
\ncline{0}{1}\rput(1,.3){\scriptsize$Q_1$}
\end{pspicture}
\begin{pspicture}(-.5,-.5)(3,3)
\cnode(1,1){.1}{00}\rput(.7,1){\scriptsize00}
\cnode(1,3){.1}{10}\rput(.7,3){\scriptsize01}
\cnode(3,1){.1}{01}\rput(3.3,1){\scriptsize10}
\cnode(3,3){.1}{11}\rput(3.3,3){\scriptsize11}
\ncline{00}{01}\ncline{01}{11}\ncline{11}{10}\ncline{10}{00}
\rput(2,.3){\scriptsize$Q_2$}
\end{pspicture}
\begin{pspicture}(-.5,-.5)(3,4)
\cnode(1,1){.1}{000}\rput(.64,1){\scriptsize000}
\cnode(1,3){.1}{001}\rput(.64,3){\scriptsize001}
\cnode(3,1){.1}{100}\rput(3.4,1){\scriptsize100}
\cnode(3,3){.1}{101}\rput(3.4,3){\scriptsize101}
\cnode(1.7,1.7){.1}{010}\rput(1.35,1.7){\scriptsize010}
\cnode(1.7,3.7){.1}{011}\rput(1.35,3.7){\scriptsize011}
\cnode(3.7,1.7){.1}{110}\rput(4.05,1.7){\scriptsize110}
\cnode(3.7,3.7){.1}{111}\rput(4.05,3.7){\scriptsize111}
\ncline{000}{001}\ncline{001}{101}\ncline{101}{100}\ncline{100}{000}
\ncline{010}{011}\ncline{011}{111}\ncline{111}{110}\ncline{110}{010}
\ncline{000}{010}\ncline{001}{011}\ncline{101}{111}\ncline{100}{110}
\rput(2,.3){\scriptsize$Q_3$}
\end{pspicture}
\vskip2pt
\begin{pspicture}(-3.,0)(8,4)
\cnode(1,1){.1}{0000}\rput(.6,1){\scriptsize0000}
\cnode(1,3){.1}{0100}\rput(.6,3){\scriptsize0100}
\cnode(3,1){.1}{0001}\rput(2.66,1.2){\scriptsize0001}
\cnode(3,3){.1}{0101}\rput(2.66,3.2){\scriptsize0101}
\cnode(1.9,1.8){.1}{0010}\rput(1.5,1.85){\scriptsize0010}
\cnode(1.9,3.8){.1}{0110}\rput(1.5,3.85){\scriptsize0110}
\cnode(3.9,1.8){.1}{0011}\rput(3.55,2){\scriptsize0011}
\cnode(3.9,3.8){.1}{0111}\rput(3.55,4){\scriptsize0111}
\cnode(5,1){.1}{1001}\rput(5.4,.8){\scriptsize1001}
\cnode(5,3){.1}{1101}\rput(5.4,2.8){\scriptsize1101}
\cnode(7,1){.1}{1000}\rput(7.4,.85){\scriptsize1000}
\cnode(7,3){.1}{1100}\rput(7.4,2.85){\scriptsize1100}
\cnode(5.9,1.8){.1}{1011}\rput(6.3,2){\scriptsize1011}
\cnode(5.9,3.8){.1}{1111}\rput(6.3,4){\scriptsize1111}
\cnode(7.9,1.8){.1}{1010}\rput(8.34,1.8){\scriptsize1010}
\cnode(7.9,3.8){.1}{1110}\rput(8.34,3.8){\scriptsize1110}
\ncline{0000}{0001}\ncline{0001}{0101}\ncline{0101}{0100}\ncline{0100}{0000}
\ncline{0010}{0011}\ncline{0011}{0111}\ncline{0111}{0110}\ncline{0110}{0010}
\ncline{0000}{0010}\ncline{0001}{0011}\ncline{0101}{0111}\ncline{0100}{0110}
\ncline{1001}{1000}\ncline{1000}{1010}\ncline{1010}{1011}\ncline{1011}{1001}
\ncline{1101}{1100}\ncline{1100}{1110}\ncline{1110}{1111}\ncline{1111}{1101}
\ncline{1001}{1101}\ncline{1000}{1100}\ncline{1010}{1110}\ncline{1011}{1111}
\nccurve[angleA=-20,angleB=-160]{0000}{1000}
\nccurve[angleA=-20,angleB=-160]{0010}{1010}
\nccurve[angleA=-20,angleB=-160]{0001}{1001}
\nccurve[angleA=-20,angleB=-160]{0011}{1011}
\nccurve[angleA=20,angleB=160]{0100}{1100}
\nccurve[angleA=20,angleB=160]{0110}{1110}
\nccurve[angleA=20,angleB=160]{0101}{1101}
\nccurve[angleA=20,angleB=160]{0111}{1111}
\rput(4,.1){\scriptsize$Q_4$}
\end{pspicture}
\caption{\label{f10}                                                           
\footnotesize {The hypercubes $Q_n$, where $Q_1=K_2$, $Q_i=K_2\times
Q_{i-1}$ for $i=2,3,4$}}
\end{figure}
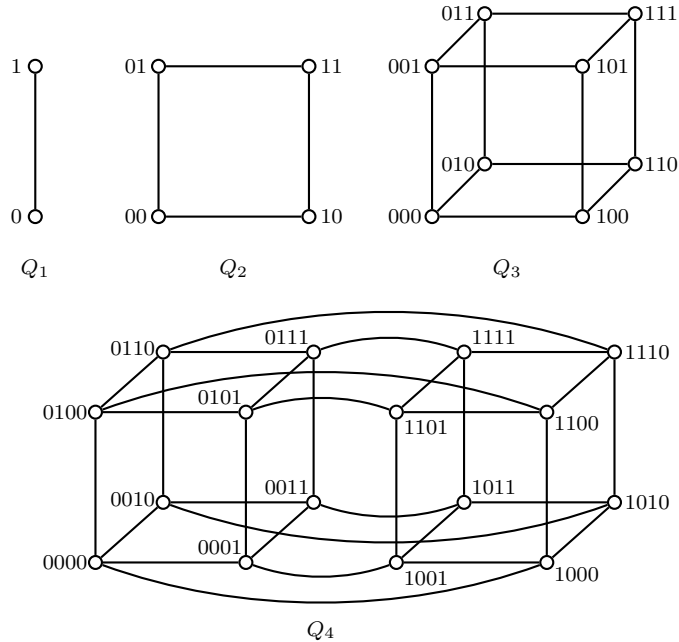

Examples of the Cartesian product are shown in Figure~\ref{f10},
where $Q_1=K_2$, $Q_i=K_2\times Q_{i-1}$ for $i=2,3,4$.


As an operation of graphs, the Cartesian products satisfy
commutative and associative laws if we identify isomorphic graphs.
Thus, we can define the Cartesian product $G_1\times G_2\times
\cdots \times  G_n$. There is an edge between a vertex $x_1x_2\cdots
x_n$ and another $y_1y_2\cdots y_n$ if and only if they differ
exactly in the $i$th coordinate and $x_iy_i\in E(G_i)$.

\vskip6pt

The Cartesian product $\Gamma=\Gamma_1\times \Gamma_2\times\cdots
\times \Gamma_n=(X, \circ )$ of $n$ finite groups $\Gamma_i=(X_i,
\circ_i)$ for each $i=1,2,\ldots,n$, where $X=X_1\times
X_2\times\cdots \times X_n$. The operation $\circ$ is defined as
follows：
$$
(x_1x_2\cdots x_n)\circ (y_1y_2\cdots y_n)
=(x_1\circ_1y_1)(x_2\circ_2y_2)\cdots(x_n\circ_ny_n),
$$
where $x_i, \ y_i\in X_i\ (i=1, 2, \ldots, n)$. For $x_1x_2\cdots
x_n\in\Gamma$, its inverse $(x_1x_2\cdots$
$x_n)^{-1}=x_1^{-1}x_2^{-1}\cdots x_n^{-1}$, the identity
$e=e_1e_2\cdots e_n$, where $x_i^{-1}$ is the inverse of $x_i$ in
$\Gamma_i$, $e_i$ is the identity in $\Gamma_i$ for each $i=1, 2,
\ldots, n$.

For example, consider $Z_4\times Z_2=\{00, 10, 20, 30, 01, 11, 21,
31\}$. For any $x_1x_2, y_1y_2\in Z_4\times Z_2$, $x_1, y_1\in Z_4,
\ x_2, y_2\in Z_2$, definite the operation:
$$
(x_1x_2)\circ (y_1y_2) =(x_1+y_1)(\mbox{mod}\
4)(x_2+y_2)(\mbox{mod}\ 2).
$$
It is easy to verify that under the above operation, $Z_4\times Z_2$
forms a group, the identity is $00$.

Consider the additive group $Z_k (k\geqslant 2)$ of residue classes
modulo $k$, that is the ring group with order $k$, zero is the
identity, the inverse of $i$ is $k-i$. If $S=\{1\}$, then $S^{-1}=S$
for $k=2$; and $S^{-1}\ne S$ otherwise. Thus the Cayley graph
$C_{Z_2}(\{1\})=K_2$, the Cayley graph $C_{Z_k}(\{1,k-1\})$ is a
cycle $C_k$ if $k\geqslant 3$.

\begin{lem}{\rm\cite{x13}}\label{lem3.16}\ The Cartesian product of
Cayley graphs is a Cayley graph. More precisely speaking, let
$G_i=C_{\Gamma _i}(S_i)$ be a Cayley graph of a finite group $\Gamma
_i$ with respect to a subset $S_i$, then $G=G_1\times G_2\times
\cdots \times  G_n$ is a Cayley graph $C_{\Gamma}(S)$ of the group
$\Gamma =\Gamma _1\times \Gamma _2 \times \cdots \times \Gamma _n$
with respect to the subset
$$
 S=\bigcup _{i=1}^n\{e_1\cdots e_{i-1}\}\times S_i\times
 \{e_{i+1}\cdots e_n\},
$$
where $e_i$ is the identity of $\Gamma _i$ for each
$i=1,2,\ldots,n$.\end{lem}

Let $\Gamma$ be the Cartesian product of $n (\geqslant 2)$ additive
groups $Z_k$, i.e., $\Gamma =Z_k\times Z_k\times \cdots \times Z_k$,
and let
$$
 S=\bigcup _{i=1}^n\{e_1\cdots e_{i-1}\}\times S_i\times
 \{e_{i+1}\cdots e_n\},
$$
where $e_i=0$ and $S_i=\{1,k-1\}$ for each $i=1,2,\ldots,n$. By
Lemma~\ref{lem3.16}, $C_{\Gamma}(S)$ is a Cayley graph. For example,
let $k=2$, then
$$
\begin{array}{rl}
 S &=\displaystyle\bigcup_{i=1}^n\{e_1\cdots e_{i-1}\}\times
                   S_i\times \{e_{i+1}\cdots e_n\}\\
 & =\{100\cdots 00, 010\cdots 00,\ldots, 000\cdots 01\},
 \end{array}
$$
where $S_i=\{1\}$ for $i=1,2,\ldots,n$. The Cayley graph
$C_{\Gamma}(S)=\underbrace{K_2\times K_2\times \cdots \times K_2}_n$
is the well-known hypercube $Q_n$.

\vskip6pt

When $k\geqslant 3$, the Cayley graph
$C_{\Gamma}(S)=\underbrace{C_{k}\times C_{k}\times \cdots \times
C_{k}}_n$ is called the {\it $k$-ary $n$-cube}, first studied by
Dally~\cite{d90} and denoted by $Q_{n}^{k}$ (also see
Xu~\cite{x13}), which is an $2n$-regular graph with $k^n$ vertices
and $n\,k^n$ edges.

\begin{lem}{\rm \cite{gh14,hlh07}}\label{lem3.17}
For any $x,y\in V(Q_{n}^{k})$, $k\geqslant 2$,
 $$
 |N(x)\cap N(y)|=\left\{\begin{array}{ll}
 1& \ \text{if $xy\in E(Q_{n}^k)$ and $k=3$};\\
 2& \ \text{if $xy\notin E(Q_{n}^k)$ and $N(x)\cap N(y)\ne\emptyset$};\\
 0& \ \text{otherwise.}
 \end{array}\right.
 $$
\end{lem}

\begin{lem}{\rm \cite{ghl13, gh14, hcsht09}}\label{lem3.18}
Let $F$ be a vertex-cut of $Q_{n}^{k}$ ($n\geqslant 5$) with
 $$
 |F|\leqslant\left\{\begin{array}{ll}
  6n-6\ & \text{if $k\geqslant 4$};\\
  6n-8\ & \text{if $k=3$}; \\
  3n-6\ & \text{if $k=2$}.
  \end{array}\right.
  $$
Then $Q_{n}^{k}-F$ has a large component and small components have
at most two vertices in total.
  \end{lem}

Xu {\it et al.}~\cite{xzhz05} determined $\kappa_2(Q^2_n)=3n-5$ for
$n\geqslant 4$. Zhao and Jin~\cite{zj13} determined
$\kappa_2(Q^3_n)=6n-7$ for $n\geqslant 3$. Hsieh {\it et
al.}~\cite{hc12} determined $\kappa_2(Q^k_n)=6n-5$ for $k\geqslant
4$ and $n\geqslant 5$. Hsu {\it et al.}~\cite{hcsht09} proved
$t_{c}(Q_{n}^{2})=3n-5$ for $n\geqslant 5$. By Theorem~\ref{thm2.4},
we immediately obtain the following result which contains the above
results.

 \begin{thm}\label{thm3.19}
For $n\geqslant 8$ if $k=5$ and $n\geqslant 6$ otherwise,
$t_{c}(Q_{n}^{k})=t=\kappa_{2}(Q_{n}^{k})$, where
$$
  t=\left\{\begin{array}{ll}
  6n-5\ & \text{if $k\geqslant 4$};\\
  6n-7\ & \text{if $k=3$}; \\
  3n-5\ & \text{if $k=2$}.
  \end{array}\right.
  $$
 \end{thm}

\begin{pf}
Note that $Q_{n}^{k}$ is $n$-regular for $k=2$, and $2n$-regular for
$k\geqslant 3$, and $Q_{n}^{k}$ contains $C_3$ if and only if $k=3$
and contains $C_5$ if and only if $k=5$. By Lemma~\ref{lem3.17}, it
is easy to verify that $t=\min\{|N(T)|:\ T=P_3\ {\rm or}\ C_3$ in
$Q_{n}^{k}\}=|N(P_3)|$, where $P_3$ is any 3-path in $Q_{n}^{k}$
since $Q_{n}^{k}$ is vertex-transitive.

Let $F=N(P_3)$. Then $F$ is a vertex-cut of $Q_{n}^{k}$ and $|F|=t$.
To prove the theorem, we only need to verify that $Q_{n}^{k}$
satisfies conditions in Theorem~\ref{thm2.4}.

{\rm (a)}\ If $|F|\leqslant t-1$, then by Lemma~\ref{lem3.18},
$Q_n^k-H$ has a large component and small components has at most two
vertices in total.

{\rm (b)}\ By Lemma~\ref{lem3.17}, $n\geqslant 3\ell(Q_n^k)+2=8$ if
$k=5$, and $n\geqslant 2\ell(Q_n^k)+2=6$ otherwise.

{\rm (c)}\ For $n\geqslant 8$, it is easy to verify that
 $$
 |V(Q_{n}^{k})|-(\Delta+1)(t-1)-4=\left\{
 \begin{array}{ll}
 2^n-(n+1)(3n-6)-4>0\ & {\rm if}\ k=2;\\
 3^n-(2n+1)(6n-8)-4>0\ & {\rm if}\ k=3;\\
 k^n-(2n+1)(6n-6)-4>0\ & {\rm if}\ k\geqslant 4.
 \end{array}\right.
 $$

Thus, $Q_n^k$ satisfies all conditions in Theorem~\ref{thm2.4}, and
so $t_c(Q_n^k)=t=\kappa_2(Q_n^k)$ for $n\geqslant 8$ if $k=5$ and
$n\geqslant 6$ otherwise.
\end{pf}

\begin{cor}
$t_{c}(Q_{n}^{2})=3n-5=\kappa_2(Q^2_n)$ and
$t_{c}(Q_{n}^{3})=6n-7=\kappa_2(Q^3_n)$ for $n\geqslant 6$
\end{cor}

\subsection{Dual-Cubes}

A dual-cube $DC_n$, proposed by Li and Peng~\cite{p00}, consists of
$2^{2n+1}$ vertices, and each vertex is labeled with a unique
($2n+1$)- bits binary string and has $n+1$ neighbors. There is a
link between two nodes $u=u_{2n}u_{2n-1}\ldots u_0$ and
$v=v_{2n}v_{2n-1}\ldots v_0$ if and only if $u$ and $v$ differ
exactly in one bit position $i$ under the the following conditions:

(a) if $0\leq i \leq n-1$, then $u_{2n}=v_{2n}=0$; and

(b) if $n\leq i \leq 2n-1$, then $u_{2n}=v_{2n}=1$.

\begin{figure}[h]
\psset{unit=.8}
\begin{pspicture}(-5,-1.4)(1.7,5.)
 \cnode(-1.5,4){.09}{11} \cnode(0,4){.09}{12}
 \cnode(-1.5,2.5){.09}{14} \cnode(0,2.5){.09}{13}
 \ncline{11}{12}\ncline{12}{13}\ncline{13}{14}\ncline{14}{11}

 \cnode(1.5,4){.09}{21} \cnode(3,4){.09}{22}
 \cnode(1.5,2.5){.09}{24} \cnode(3,2.5){.09}{23}
 \ncline{21}{22}\ncline{22}{23}\ncline{23}{24}\ncline{24}{21}

 \cnode(4.5,4){.09}{31} \cnode(6,4){.09}{32}
 \cnode(4.5,2.5){.09}{34} \cnode(6,2.5){.09}{33}
 \ncline{31}{32}\ncline{32}{33}\ncline{33}{34}\ncline{34}{31}

 \cnode(7.5,4){.09}{41} \cnode(9,4){.09}{42}
 \cnode(7.5,2.5){.09}{44} \cnode(9,2.5){.09}{43}
 \ncline{41}{42}\ncline{42}{43}\ncline{43}{44}\ncline{44}{41}

 \cnode(-1.5,1){.09}{51} \cnode(0,1){.09}{52}
 \cnode(-1.5,-.5){.09}{54} \cnode(0,-.5){.09}{53}
 \ncline{51}{52}\ncline{52}{53}\ncline{53}{54}\ncline{54}{51}

 \cnode(1.5,1){.09}{61} \cnode(3,1){.09}{62}
 \cnode(1.5,-.5){.09}{64} \cnode(3,-.5){.09}{63}
 \ncline{61}{62}\ncline{62}{63}\ncline{63}{64}\ncline{64}{61}

 \cnode(4.5,1){.09}{71} \cnode(6,1){.09}{72}
 \cnode(4.5,-.5){.09}{74} \cnode(6,-.5){.09}{73}
 \ncline{71}{72}\ncline{72}{73}\ncline{73}{74}\ncline{74}{71}

 \cnode(7.5,1){.09}{81} \cnode(9,1){.09}{82}
 \cnode(7.5,-.5){.09}{84} \cnode(9,-.5){.09}{83}
 \ncline{81}{82}\ncline{82}{83}\ncline{83}{84}\ncline{84}{81}

 \ncline{14}{61}\ncline{13}{71}\ncline{24}{51}\ncline{23}{81}
 \ncline{34}{52}\ncline{33}{82}\ncline{44}{62}\ncline{43}{72}

 \psarc(-0.05,0.5){1.2in}{68}{112}\psarc(7.45,0.5){1.2in}{68}{112}
 \psarc(2.2,-1.25){1.8in}{67}{112}\psarc(5.2,-1.25){1.8in}{67}{112}

 \psarc(-0.05,3.){1.2in}{249}{293} \psarc(7.45,3.){1.2in}{249}{293}
 \psarc(2.2,4.7){1.8in}{248}{293}\psarc(5.2,4.7){1.8in}{248}{293}

 \rput(-2.1,4.){\scriptsize $00000$}\rput(-2.1,2.5){\scriptsize $00010$}
 \rput(-0.55,3.75){\scriptsize $00001$}\rput(-0.55,2.75){\scriptsize $00011$}
 \rput(-2.1,1.){\scriptsize $01000$}\rput(-2.1,-0.5){\scriptsize $01010$}
 \rput(-0.55,.8){\scriptsize $01001$}\rput(-0.55,-0.25){\scriptsize $01011$}

 \rput(.95,3.75){\scriptsize $10000$}\rput(.95,2.75){\scriptsize $11000$}
 \rput(.95,.8){\scriptsize $10010$}\rput(.95,-0.25){\scriptsize $11010$}
 \rput(2.45,3.75){\scriptsize $10100$}\rput(2.45,2.75){\scriptsize $11100$}
 \rput(2.45,.8){\scriptsize $10110$}\rput(2.45,-0.25){\scriptsize $11110$}

 \rput(3.95,3.75){\scriptsize $10001$}\rput(3.95,2.75){\scriptsize $11001$}
 \rput(3.95,.8){\scriptsize $10011$}\rput(3.95,-0.25){\scriptsize $11011$}
 \rput(5.45,3.75){\scriptsize $10101$}\rput(5.45,2.75){\scriptsize $11101$}
 \rput(5.45,.8){\scriptsize $10111$}\rput(5.45,-0.25){\scriptsize $11111$}

 \rput(6.95,3.75){\scriptsize $00100$}\rput(6.95,2.75){\scriptsize $00110$}
 \rput(6.95,.8){\scriptsize $01100$}\rput(6.95,-0.25){\scriptsize $01110$}
 \rput(9.65,4.1){\scriptsize $00101$}\rput(9.65,2.5){\scriptsize $00111$}
 \rput(9.65,1.1){\scriptsize $01101$}\rput(9.65,-0.4){\scriptsize $01111$}

\end{pspicture}
\caption{ \label{f11} The dual-cube $DC_2$}
\end{figure}
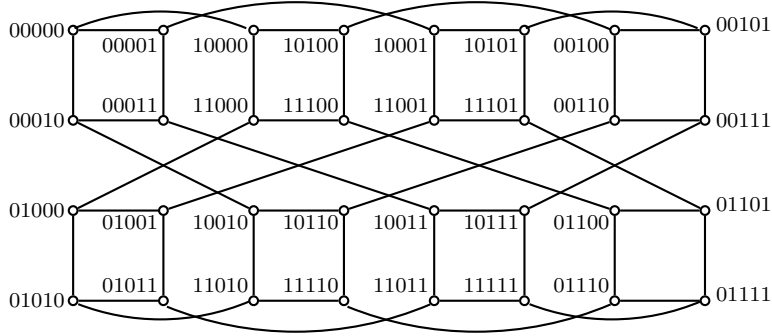

Figure~\ref{f11} shows the bubble-sort graphs $DC_2$. A dual-cube
$DC_n$ is an $(n+1)$-regular bipartite graph of order $2^{2n+1}$.
Moreover,
%
%
Zhou {\it et al.}~\cite{zcx12} showed that $DC_n$ is a Cayley graph,
and so $DC_n$ is vertex-transitive.

\begin{lem}\label{lem3.21}\textnormal{(Zhou {\it et al.}~\cite{zcx12})}
For any $x,y\in V(DC_n)$, if $xy\notin E(DC_n)$ and $N(x)\cap
N(y)\ne\emptyset$, then $|N(x)\cap N(y)|\leqslant 2$.
\end{lem}

Since $DC_n$ is an $(n+1)$-regular bipartite graph, and so it
contains no $C_3$, according to Lemma~\ref{lem3.21}, if
$P_3=(x,y,z)$ is a 3-path, where $xz\notin E(G)$, then $|N(x)\cap
N(y)|=|N(y)\cap N(z)|=0$ and $|N(x)\cap N(z)|\leqslant 2$, and so
the number of neighbors of $P_3$ in $DC_n$ can be counted as
follows.
 $$
 \begin{array}{rl}
 |N(P_3)|&=d(x)+d(y)+d(z)-|N(x)\cap N(z)|-\Sigma(P_3)\\
  &=3(n+1)-|N(x)\cap N(z)|-4\\
  &=\left\{\begin{array}{ll}
  3n-1 & \ {\rm if} |N(x)\cap N(z)|=1;\\
  3n-2 & \ {\rm if} |N(x)\cap N(z)|=2.
  \end{array}\right.
  \end{array}
 $$
Since $DC_n$ is vertex-transitive, for any $3$-path $P_3$ in $DC_n$
that $|N(P_3)|$ is the smallest, we have that
  \begin{equation}\label{e6}
  |N(P_3)|=3n-2.
   \end{equation}

\begin{lem}\label{lem3.22} \textnormal{(Zhou {\it et al.}~\cite{zcx12})}
Let $F\subset V(DC_n)$ with $|F|\leqslant 3n-3$ and $n\geqslant 3$.
If $DC_n-F$ is disconnected, then it has either two components, one
of which is an isolated vertex or an edge, or three components, two
of which are isolated vertices.
\end{lem}

Zhou {\it et al.}~\cite{zcx12} determined $\kappa_2(DC_n)=3n-2$ and
$t_c(DC_n)=3n-2$ for $n\geq 3$, dependently. By
Theorem~\ref{thm2.4}, we immediately obtain the following result
which contains the above results.


\begin{thm}\label{thm3.23}\
$t_c(DC_n)=3n-2=\kappa_2(DC_n)$ for $n\geqslant 5$.
\end{thm}

\begin{pf}
Since $DC_n$ contains no $C_3$, $t=\min\{|N(T)|:\ T=P_3$ or $C_3$ in
$DC_n\}=|N(P_3)|$, where $P_3$ is any 3-path in $DC_n$ since $DC_n$
is vertex-transitive. Let $F=N(P_3)$. Then $|F|=t=3n-2$ by
(\ref{e6}). It is easy to check that $F$ is a vertex-cut of $DC_n$.
To prove the theorem, we only need to verify that $DC_n$ satisfies
conditions in Theorem~\ref{thm2.4}.

{\rm (a)}\ If $|F|\leqslant t-1$ then, by Lemma~\ref{lem3.22},
$DC_n-F$ has a large component and small components which contain at
most two vertices in total.

{\rm (b)}\ By Lemma~\ref{lem3.21}, $\ell(DC_n)=2$. Since $DC_{n}$ is
$(n+1)$-regular bipartite, it contains no $5$-cycle, and so
$n+1\geqslant 6=2\ell(DC_n)+2$.

{\rm (c)}\ It is easy to check that
 $2^{2n+1}-(n+2)(t-1)-4=2^{2n+1}-(n+2)(3n-3)-4>0$ for $n\geqslant 5$.

$DC_n$ satisfies all conditions in Theorem~\ref{thm2.4}, and so
$t_c(DC_n)=3n-2=\kappa_2(DC_n)$. The theorem follows.
\end{pf}

\section{Conclusions}

The conditional diagnosability $t_c(G)$ under the comparison model
and the 2-extra connectivity $\kappa_2(G)$ are two important
parameters to measure ability of diagnosing faulty processors and
fault-tolerance in a multiprocessor system $G$ with the presence of
failing processors. Although these two parameters have attracted
considerable attention and determined for many classes of well-known
graphs in recent years, but are obtained independently. This paper
establishes the close relationship between these two parameters by
proving $t_c(G)=\kappa_2(G)$ for a regular graph $G$ with some
acceptable conditions. As applications, the conditional
diagnosability and the 2-extra connectivity are determined for some
well-known classes of vertex-transitive graphs such as star graphs,
$(n,k)$-star graphs, $(n,k)$-arrangement graphs, Cayley graphs
obtained from transposition generating trees, $k$-ary $n$-cube
networks and dual-cubes. Furthermore, many known results about these
networks are obtained directly.

Under the comparison diagnosis model, the diagnosability and the
1-extra connectivity should have some relationships. On the other
hand, in addition to the comparison diagnosis model, there are
several other diagnosis models such as the PMC model. Under the PMC
model, what is the relationship between the diagnosability or the
conditional diagnosability and the $h$-extra connectivity for some
$h$? These will be explored in future.

\end{document}